\documentclass[11pt,twoside]{article}
\usepackage{amsmath, amsthm, amscd, amsfonts, amssymb, graphicx, color}

\date{}

\begin{document}

\centerline{}

\centerline {\Large{\bf Introduction to Continuous biframes in Hilbert spaces}}
\centerline {\Large{\bf and their tensor products}}

\newcommand{\mvec}[1]{\mbox{\bfseries\itshape #1}}
\centerline{}
\centerline{\textbf{Prasenjit Ghosh}}
\centerline{Department of Mathematics,}
\centerline{Barwan N. S. High School (HS), Barwan, }
\centerline{Murshidabad, 742161, West Bengal, India}
\centerline{e-mail: prasenjitpuremath@gmail.com}
\centerline{}
\centerline{\textbf{T. K. Samanta}}
\centerline{Department of Mathematics, Uluberia College,}
\centerline{Uluberia, Howrah, 711315,  West Bengal, India}
\centerline{e-mail: mumpu$_{-}$tapas5@yahoo.co.in}

\newtheorem{Theorem}{\quad Theorem}[section]

\newtheorem{definition}[Theorem]{\quad Definition}

\newtheorem{theorem}[Theorem]{\quad Theorem}

\newtheorem{remark}[Theorem]{\quad Remark}

\newtheorem{corollary}[Theorem]{\quad Corollary}

\newtheorem{note}[Theorem]{\quad Note}

\newtheorem{lemma}[Theorem]{\quad Lemma}

\newtheorem{example}[Theorem]{\quad Example}

\newtheorem{result}[Theorem]{\quad Result}
\newtheorem{conclusion}[Theorem]{\quad Conclusion}

\newtheorem{proposition}[Theorem]{\quad Proposition}

\begin{abstract}
\textbf{\emph{We introduce the notion of a continuous biframe in a Hilbert space which is a generalization of discrete biframe in Hilbert space. Representation theorem for this type of generalized frame is verified and some characterizations of this biframe with the help of a invertible operator is given.\,Here we also introduce the concept of continuous biframe for the tensor products of Hilbert spaces and give an example.\,Further, we study dual continuous biframe and continuous biframe Bessel multiplier in Hilbert spaces and their tensor products. }}
\end{abstract}
\section{Introduction}
 
\smallskip\hspace{.6 cm}The notion of a frame in Hilbert space was first introduced by Duffin and Schaeffer \cite{Duffin} in connection with some fundamental problem in non-harmonic analysis.\,Thereafter, it was further developed and popularized by Daubechies et al.\,\cite{Daubechies} in 1986.\,A discrete frame is a countable family of elements in a separable Hilbert space which allows for a stable, not necessarily unique, decomposition of an arbitrary element into an expansion of the frame element.\,A sequence \,$\left\{\,f_{\,i}\,\right\}_{i \,=\, 1}^{\infty}$\, in a separable Hilbert space \,$H$\, is called a frame for \,$H$, if there exist positive constants \,$0 \,<\, A \,\leq\, B \,<\, \infty$\, such that
\[ A\; \|\,f\,\|^{\,2} \,\leq\, \sum\limits_{i \,=\, 1}^{\infty}\, \left|\ \left <\,f \,,\, f_{\,i} \, \right >\,\right|^{\,2} \,\leq\, B \,\|\,f\,\|^{\,2}\; \;\text{for all}\; \;f \,\in\, H.\]
The constants \,$A$\, and \,$B$\, are called lower and upper bounds, respectively.

Controlled frame is one of the newest generalization of frame in Hilbert space.\,I. Bogdanova et al.\,\cite{I} introduced controlled frame for spherical wavelets to get numerically more efficient approximation algorithm.\,Thereafter, P. Balaz \cite{B} developed weighted and controlled frame in Hilbert space.\,Biframe is also a generalization of controlled frame in Hilbert space which was studied by M. F. Parizi et al.\,\cite{MF}.\,To define frame in Hilbert space, only one sequence is needed, but for a biframe, two sequences are needed.\,A pair of sequences \,$\left(\,\left\{\,f_{i}\,\right\}_{i \,=\, 1}^{\,\infty}\,,\, \left\{\,g_{i}\,\right\}_{i \,=\, 1}^{\,\infty}\,\right)$\, in \,$H$\, is called a biframe for \,$H$\, if there exist positive constants \,$A$\, and \,$B$\, such that
\[A\; \|\,f\,\|^{\,2} \,\leq\, \sum\limits_{i \,=\, 1}^{\infty}\, \left <\,f\,,\, f_{\,i} \, \right >\,\left<\,g_{\,i}\,,\, f\,\right> \,\leq\, B \,\|\,f\,\|^{\,2}\; \;\forall\; f \,\in\, H.\] 
The constants \,$A$\, and \,$B$\, are called lower and upper biframe bounds, respectively.

\,Continuous frames extend the concept of discrete frames when the indices are related to some measurable space.\,Continuous frame in Hilbert space was studied by A. Rahimi et al \cite{AR}.\,M. H. faroughi and E. Osgooei \cite{MH} also studied continuous frame and continuous Bessel mapping.\,Continuous frame and discrete frame have been used in image processing, coding theory, wavelet analysis, signal denoising, feature extraction, robust signal processing etc.

In this paper, we give the notion of continuous biframes in Hilbert spaces and their tensor products and then discuss some examples of this type of frame.\,A characterization of continuous biframe using its frame operator is established.\,We will see that the image of a continuous biframe under a bounded invertible operator in Hilbert space is also a continuous biframe in Hilbert space.\,Continuous biframe Bessel multipliers in Hilbert spaces and their tensor product are presented. 
 
\section{Preliminaries}

\begin{definition}\cite{AR}
Let \,$H$\, be a complex Hilbert space and \,$(\,\Omega,\, \mu\,)$\, be a measure space with positive measure \,$\mu$.\,A mapping \,$F \,:\, \Omega \,\to\, H$\, is called a continuous frame with respect to \,$\left(\,\Omega,\, \mu\,\right)$\, if
\begin{itemize}
\item[$(i)$] \,$F$\, is weakly-measurable, i.\,e., for all \,$f \,\in\, H$, \,$w \,\to\, \left<\,f,\, F\,(\,w\,)\,\right>$\, is a measurable function on \,$\Omega$,
\item[$(ii)$]there exist constants \,$0 \,<\, A \,\leq\, B \,<\, \infty$\, such that
\end{itemize}
\[A\,\left\|\,f\,\right\|^{\,2} \leq \int\limits_{\Omega}\left|\,\left<\,f,\, F\,(\,w\,)\,\right>\,\right|^{\,2}\,d\mu \leq B\left\|\,f\,\right\|^{\,2}\,,\]
for all \,$f \,\in\, H$.\,The constants \,$A$\, and \,$B$\, are called continuous frame bounds.\,If \,$A \,=\, B$, then it is called a tight continuous frame.\,If the mapping \,$F$\, satisfies only the right inequality, then it is called continuous Bessel mapping with Bessel bound \,$B$.
\end{definition}

Let \,$L^{\,2}\,(\,\Omega,\,\mu\,)$\, be the class of all measurable functions \,$f \,:\, \Omega \,\to\, H$\, such that \,$\|\,f\,\|^{\,2}_{\,2} \,=\,  \int\limits_{\,\Omega}\,\left\|\,f\,(\,w\,)\,\right\|^{\,2}\,d\mu \,<\, \infty$.\,It can be proved that \,$L^{\,2}\,(\,\Omega,\,\mu\,)$\, is a Hilbert space with respect to the inner product defined by
\[\left<\,f,\, g\right>_{L^{2}} \,=\, \int\limits_{\Omega}\,\left<\,f\,(\,w\,),\, g\,(\,w\,)\right>\,d\mu\,, \,f,\, g \,\in\, L^{\,2}\,(\,\Omega,\,\mu\,).\]

\begin{theorem}\cite{MH}
Let \,$F \,:\, \Omega \,\to\, H$\, be a Bessel mapping.\;Then the operator \,$T_{C} \,:\, L^{2}\left(\,\Omega,\,\mu\,\right) \,\to\, H$\, is defined by
\[\left<\,T_{C}\,(\,\varphi\,),\, h\,\right> \,=\, \int\limits_{\,\Omega}\,\varphi\,(\,w\,)\,\left<\,F\,(\,w\,),\, h\,\right>\,d\mu\]where \,$\varphi \,\in\, L^{2}\left(\,\Omega,\,\mu\,\right)$\, and \,$h \,\in\, H$\, is well-defined, linear, bounded and its adjoint operator is given by  
\[T^{\,\ast}_{C} \,:\, H \,\to\, L^{2}\left(\,\Omega,\,\mu\,\right) \;,\; T^{\,\ast}_{C}\,f\,(\,w\,) \,=\, \left<\,f,\, F\,(\,w\,)\,\right>\;,\; f \,\in\, H\,,\;\; w \,\in\, \Omega.\]
\end{theorem}

The operator \,$T_{C}$\, is called a pre-frame operator or synthesis operator and its adjoint operator is called analysis operator of \,$F$.\,

\begin{definition}\cite{MH}
Let \,$F \,:\, \Omega \,\to\, H$\, be a continuous frame for \,$H$.\,Then the operator \,$S_{C} \,:\, H \,\to\, H$\, defined by
\[\left<\,S_{C}\,(\,f\,),\, h\,\right> \,=\, \int\limits_{\,\Omega}\,\left<\,f, F\,(\,w\,)\,\right>\left<\,F\,(\,w\,),\, h\,\right>\,d\mu\,, \;\forall\, \,f,\, h \,\in\, H\]is called the frame operator of \,$F$.
\end{definition}

The tensor product of Hilbert spaces are introduced by several ways and it is a certain linear space of operators which was represented by Folland in \cite{Folland}.

\begin{definition}\cite{Upender}\label{1.def1.01}
The tensor product of Hilbert spaces \,$H_{1}$\, and \,$H_{2}$\, is denoted by \,$H_{1} \,\otimes\, H_{2}$\, and it is defined to be an inner product space associated with the inner product
\begin{equation}\label{eq1.001}   
\left<\,f \,\otimes\, g \,,\, f^{\,\prime} \,\otimes\, g^{\,\prime}\,\right> \,=\, \left<\,f,\, f^{\,\prime}\,\right>_{\,1}\;\left<\,g,\, g^{\,\prime}\,\right>_{\,2},
\end{equation}
for all \,$f,\, f^{\,\prime} \,\in\, H_{1}\; \;\text{and}\; \;g,\, g^{\,\prime} \,\in\, H_{2}$.\,The norm on \,$H_{1} \,\otimes\, H_{2}$\, is given by 
\begin{equation}\label{eq1.0001}
\left\|\,f \,\otimes\, g\,\right\| \,=\, \|\,f\,\|_{\,1}\;\|\,g\,\|_{\,2}\; \;\forall\; f \,\in\, H_{1}\; \;\text{and}\; \,g \,\in\, H_{2}.
\end{equation}
The space \,$H_{1} \,\otimes\, H_{2}$\, is complete with respect to the above inner product.\;Therefore the space \,$H_{1} \,\otimes\, H_{2}$\, is a Hilbert space.     
\end{definition} 

For \,$Q \,\in\, \mathcal{B}\,(\,H_{1}\,)$\, and \,$T \,\in\, \mathcal{B}\,(\,H_{2}\,)$, the tensor product of operators \,$Q$\, and \,$T$\, is denoted by \,$Q \,\otimes\, T$\, and defined as 
\[\left(\,Q \,\otimes\, T\,\right)\,A \,=\, Q\,A\,T^{\,\ast}\; \;\forall\; \;A \,\in\, H_{1} \,\otimes\, H_{2}.\]
It can be easily verified that \,$Q \,\otimes\, T \,\in\, \mathcal{B}\,(\,H_{1} \,\otimes\, H_{2}\,)$\, \cite{Folland}.\\

\begin{theorem}\cite{Folland}\label{th1.1}
Suppose \,$Q,\, Q^{\prime} \,\in\, \mathcal{B}\,(\,H_{1}\,)$\, and \,$T,\, T^{\prime} \,\in\, \mathcal{B}\,(\,H_{2}\,)$, then \begin{itemize}
\item[$(i)$] \,$Q \,\otimes\, T \,\in\, \mathcal{B}\,(\,H_{1} \,\otimes\, H_{2}\,)$\, and \,$\left\|\,Q \,\otimes\, T\,\right\| \,=\, \|\,Q\,\|\; \|\,T\,\|$.
\item[$(ii)$] \,$\left(\,Q \,\otimes\, T\,\right)\,(\,f \,\otimes\, g\,) \,=\, Q\,(\,f\,) \,\otimes\, T\,(\,g\,)$\, for all \,$f \,\in\, H_{1},\, g \,\in\, H_{2}$.
\item[$(iii)$] $\left(\,Q \,\otimes\, T\,\right)\,\left(\,Q^{\,\prime} \,\otimes\, T^{\,\prime}\,\right) \,=\, (\,Q\,Q^{\,\prime}\,) \,\otimes\, (\,T\,T^{\,\prime}\,)$. 
\item[$(iv)$] \,$Q \,\otimes\, T$\, is invertible if and only if \,$Q$\, and \,$T$\, are invertible, in which case \,$\left(\,Q \,\otimes\, T\,\right)^{\,-\, 1} \,=\, \left(\,Q^{\,-\, 1} \,\otimes\, T^{\,-\, 1}\,\right)$.
\item[$(v)$] \,$\left(\,Q \,\otimes\, T\,\right)^{\,\ast} \,=\, \left(\,Q^{\,\ast} \,\otimes\, T^{\,\ast}\,\right)$.  
\end{itemize}
\end{theorem}

\section{Continuous biframe in Hilbert space}
In this section, first we give the definition of a continuous biframe in Hilbert space and then discuss some of its properties.

\begin{definition}\label{def1.01}
Let \,$H$\, be a Hilbert space and \,$(\,\Omega,\, \mu\,)$\, be a measure space with positive measure \,$\mu$.\,A pair \,$(\,F,\, G\,) \,=\, \,\left(\,F \,:\, \Omega \,\to\, H,\, \,G \,:\, \Omega \,\to\, H\,\right)$\, of mappings is called a continuous biframe for \,$H$\, with respect to \,$\left(\,\Omega,\, \mu\,\right)$\, if
\begin{itemize}
\item[$(i)$] \,$F,\, G$\, are weakly-measurable, i.\,e., for all \,$f \,\in\, H$, \,$w \,\mapsto\, \left<\,f,\, F\,(\,w\,)\,\right>$\, and \,$w \,\mapsto\, \left<\,f,\, G\,(\,w\,)\,\right>$\, are measurable functions on \,$\Omega$,
\item[$(ii)$]there exist constants \,$0 \,<\, A \,\leq\, B \,<\, \infty$\, such that
\end{itemize}
\begin{align}
A\,\left\|\,f\,\right\|^{\,2} \leq \int\limits_{\Omega}\,\left<\,f,\, F\,(\,w\,)\,\right>\,\left<\,G\,(\,w\,),\, f\,\right>\,d\mu \leq\, B\left\|\,f\,\right\|^{\,2}\,,\label{3.eqq3.11}
\end{align}
for all \,$f \,\in\, H$.\,The constants \,$A$\, and \,$B$\, are called continuous biframe bounds.\,If \,$A \,=\, B$, then it is called a tight continuous biframe and it is called Parseval continuous biframe if \,$A \,=\, B \,=\, 1$.\,If the pair \,$(\,F,\, G\,)$\, satisfies only the right inequality, then it is called continuous biframe Bessel mapping with Bessel bound \,$B$. 
\end{definition}

In particular, if \,$\mu$\, is a counting measure and \,$\Omega \,=\, \mathbb{N}$, then \,$(\,F,\, G\,)$\, is called a discrete biframe for \,$H$.  

\begin{remark}
Let \,$F \,:\, \Omega \,\to\, H$\, be a mapping.\,Then according to the Definition \ref{def1.01}, we say that
\begin{itemize}
\item[$(i)$]If \,$(\,F,\, F\,)$\, is a continuous biframe for \,$H$, then \,$F$\, is a continuous frame for \,$H$.
\item[$(ii)$]If \,$U \,\in\, G\,\mathcal{B}(\,H\,)$, \,$(\,F,\, UF\,)$\, is a continuous biframe for \,$H$, then \,$F$\, is a \,$U$--controlled continuous frame for \,$H$, where $\mathcal{G}\,\mathcal{B}\,(\,H\,)$\, denotes the set of all bounded linear operators which have bounded inverse.  
\item[$(iii)$]If \,$T,\, U \,\in\, G\,\mathcal{B}(\,H\,)$, \,$(\,TF,\, UF\,)$\, is a continuous biframe for \,$H$, then \,$F$\, is a \,$(\,T,\, U\,)$--controlled continuous frame for \,$H$.  
\end{itemize} 
\end{remark}

Now, we validates the above definition by some examples.

\begin{example}
Let \,$H \,=\, \mathbb{R}^{\,3}$\, and \,$\left\{\,e_{\,1},\,e_{\,2},\, e_{\,3}\,\right\}$\, be an standard orthonormal basis for \,$H$.\,Consider
\[\Omega \,=\, \left\{\,x \,\in\, \mathbb{R}^{\,3} \,:\, \|\,x\,\| \,\leq\, 1\,\right\}.\]
Then it is a measure space equipped with the Lebesgue measure \,$\mu$.\,Suppose \,$\left\{\,B_{\,1},\,B_{\,2},\, B_{\,3}\,\right\}$\, is a partition of \,$\Omega$\, where \,$\mu\,(\,B_{1}\,) \,\geq\, \mu\,(\,B_{2}\,) \,\geq\, \mu\,(\,B_{3}\,) \,>\, 1$.\,Define 
\[F \,:\, \Omega \,\to\, H\hspace{.5cm}\text{by} \hspace{.3cm} F\,(\,w\,) \,=\, \begin{cases}
\dfrac{e_{\,1}}{\sqrt{\,\mu\left(\,B_{1}\,\right)}} & \text{if\;\;}\; w \,\in\, B_{1} \\ \dfrac{e_{\,2}}{\sqrt{\,\mu\left(\,B_{2}\,\right)}} & \text{if\;\;}\; w \,\in\, B_{2}\\ \dfrac{2\,e_{\,3}}{\sqrt{\,\mu\left(\,B_{3}\,\right)}} & \text{if\;\;}\; w \,\in\, B_{3} \end{cases}\]and
\[G \,:\, \Omega \,\to\, H\hspace{.5cm}\text{by} \hspace{.3cm} G\,(\,w\,) \,=\, \begin{cases}
\dfrac{2\,e_{\,1}}{\sqrt{\,\mu\left(\,B_{1}\,\right)}} & \text{if\;\;}\; w \,\in\, B_{1} \\ \dfrac{\, e_{\,2}}{\sqrt{\,\mu\left(\,B_{2}\,\right)}} & \text{if\;\;}\; w \,\in\, B_{2}\\ \dfrac{e_{\,3}}{2\,\sqrt{\,\mu\left(\,B_{3}\,\right)}} & \text{if\;\;}\; w \,\in\, B_{3} \end{cases}\]
It is easy to verify that for all \,$f \,\in\, H$, \,$w \,\mapsto\, \left<\,f,\, F\,(\,w\,)\,\right>$\, and \,$w \,\mapsto\, \left<\,f,\, G\,(\,w\,)\,\right>$\, are measurable functions on \,$\Omega$.
Now, for \,$f \,\in\, H$, we have
\begin{align*}
\int\limits_{\,\Omega}\left|\left<\,f,\, F\,(\,w\,)\,\right>\,\right|^{2}\,d\mu &= \int\limits_{\,B_{1}}\left|\,\left<\,f,\, \dfrac{1}{\sqrt{\mu\left(\,B_{1}\,\right)}}\,e_{\,1}\,\right>\,\right|^{\,2}\,d\mu \,+\, \int\limits_{\,B_{2}}\left|\,\left<\,f,\, \dfrac{1}{\sqrt{\mu\left(\,B_{2}\,\right)}}\,e_{\,2}\,\right>\,\right|^{\,2}\,d\mu\\
&\hspace{1cm}+\,\int\limits_{\,B_{3}}\left|\,\left<\,f,\, \dfrac{1}{\sqrt{\mu\left(\,B_{3}\,\right)}}\,2\,e_{\,1}\,\right>\,\right|^{\,2}\,d\mu\\
&=\,\left|\,\left<\,f,\, \,e_{\,1}\,\right>\,\right|^{\,2} \,+\, \left|\,\left<\,f,\, \,e_{\,2}\,\right>\,\right|^{\,2} \,+\, 4\,\left|\,\left<\,f,\, \,e_{\,3}\,\right>\,\right|^{\,2}\\
& \,=\, \|\,f\,\|^{\,2} \,+\, 3\,\left|\,\left<\,f,\, \,e_{\,3}\,\right>\,\right|^{\,2}.
\end{align*} 
Therefore, \,$F \,:\, \Omega \,\to\, H$\, is a continuous frame for \,$H$\, with bounds \,$1$\, and \,$4$.\,Similarly, it can be shown that \,$G \,:\, \Omega \,\to\, H$\, is a continuous frame for \,$H$\, with bounds \,$1 \,/\, 4$\, and \,$4$.\,On the other hand, for \,$f \,\in\, H$, we have
\begin{align*}
&\int\limits_{\,\Omega}\,\left<\,f,\, F\,(\,w\,)\,\right>\,\left<\,G\,(\,w\,),\, f\,\right>\,d\mu \\
&= \int\limits_{\,B_{1}}\left<\,f,\, \dfrac{e_{\,1}}{\sqrt{\mu\left(\,B_{1}\,\right)}}\,\right>\,\left<\,\dfrac{2\,e_{\,1}}{\sqrt{\mu\left(\,B_{1}\,\right)}},\, f\,\right>\,d\mu\\
& \,+\, \int\limits_{\,B_{2}}\left<\,f,\, \dfrac{e_{\,2}\,}{\sqrt{\mu\left(\,B_{2}\,\right)}}\,\right>\,\left<\,\dfrac{\, e_{\,2}}{\sqrt{\mu\left(\,B_{2}\,\right)}},\, f\,\right>\,d\mu\\
&+\,\int\limits_{\,B_{3}}\left<\,f,\, \dfrac{2\,e_{\,3}\,}{\sqrt{\mu\left(\,B_{3}\,\right)}}\,\right>\,\left<\, \dfrac{e_{\,3}\,}{2\,\sqrt{\mu\left(\,B_{3}\,\right)}},\, f\,\right>\,d\mu\\
&=\,2\,\left|\,\left<\,f,\, \,e_{\,1}\,\right>\,\right|^{\,2} \,+\, \left|\,\left<\,f,\, \,e_{\,2}\,\right>\,\right|^{\,2} \,+\, \left|\,\left<\,f,\, \,e_{\,3}\,\right>\,\right|^{\,2}\\
& \,=\, \|\,f\,\|^{\,2} \,+\, \left|\,\left<\,f,\, \,e_{\,1}\,\right>\,\right|^{\,2}.
\end{align*} 
Therefore, \,$(\,F,\, G\,)$\, is a continuous biframe for \,$H$\, with bounds \,$1$\, and \,$2$. 
\end{example}

\begin{example}\label{3exm3.11}
Let \,$H$\, be an infinite dimensional separable Hilbert space and \,$\left\{\,e_{\,i}\,\right\}_{i \,=\, 1}^{\,\infty}$\, be an orthonormal basis for \,$H$.\,Suppose
\begin{align*}
&\left\{\,f_{\,i}\,\right\}_{i \,=\, 1}^{\,\infty} \,=\, \left\{\,e_{\,1},\, e_{\,1},\, e_{\,1},\, e_{\,2},\, e_{\,3},\, \cdots\,\cdots\,\right\}\,,\\
&\left\{\,g_{\,i}\,\right\}_{i \,=\, 1}^{\,\infty} \,=\, \left\{\,0,\, e_{\,1},\, e_{\,1},\, e_{\,2},\, e_{\,3},\, \cdots\,\cdots\,\right\}.
\end{align*}
Now, let \,$\left(\,\Omega,\, \mu\,\right)$\, be a measure space with \,$\mu$\, is \,$\sigma$-finite.\,Then we can write \,$\Omega \,=\, \bigcup_{i \,=\, 1}^{\,\infty}\,\Omega_{i}$, where \,$\left\{\,\Omega_{\,i}\,\right\}_{i \,=\, 1}^{\,\infty}$\, is a sequence of disjoint measurable subsets of \,$\Omega$\, with \,$\mu\left(\,\Omega_{i}\,\right) \,<\, \infty$.\,For each \,$w \,\in\, \Omega$, we define the mappings \,$F \,:\, \Omega \,\to\, H$\, by \,$F\,(\,w\,) \,=\, \dfrac{1}{\sqrt{\mu\left(\,\Omega_{i}\,\right)}}\,f_{\,i}$\, and \,$G \,:\, \Omega \,\to\, H$\, by \,$G\,(\,w\,) \,=\, \dfrac{1}{\sqrt{\mu\left(\,\Omega_{i}\,\right)}}\,g_{\,i}$.\,Then for \,$f \,\in\, H$, we have
\begin{align*}
&\int\limits_{\,\Omega}\left|\left<\,f,\, F\,(\,w\,)\,\right>\,\right|^{2}\,d\mu \,=\,\sum\limits_{i \,=\, 1}^{\infty}\int\limits_{\,\Omega_{i}}\left|\,\left <\,f,\, f_{\,i}\,\right >\,\right|^{\,2}\,d\mu\\
& \,=\, 2\,\left|\,\left <\,f,\, e_{\,1}\,\right >\,\right|^{\,2} \,+\, \sum\limits_{i \,=\, 1}^{\infty}\, \left|\,\left <\,f,\, e_{\,i} \,\right >\,\right|^{\,2}\\
&=\,2\,\left|\,\left <\,f,\, e_{\,1}\,\right >\,\right|^{\,2} \,+\, \left\|\,f\,\right\|^{\,2}.
\end{align*} 
Therefore, \,$F$\, is a continuous frame for \,$H$\, with bounds \,$1$\, and \,$3$.\,Similarly, it can be shown that \,$G$\, is a continuous frame for \,$H$\, with bounds \,$1$\, and \,$2$.\,Now, for \,$f \,\in\, H$, we have
\begin{align*}
&\int\limits_{\,\Omega}\,\left<\,f,\, F\,(\,w\,)\,\right>\,\left<\,G\,(\,w\,),\, f\,\right>\,d\mu \,=\,\sum\limits_{i \,=\, 1}^{\infty}\int\limits_{\,\Omega_{i}}\,\left <\,f,\, f_{\,i}\,\right >\,\left <\,g_{\,i},\, f\,\right >\,d\mu\\
&=\,\left <\,f,\, e_{\,1}\,\right >\,\left <\,e_{\,1},\, f\,\right > \,+\, \left <\,f,\, e_{\,1}\,\right >\,\left <\,e_{\,1},\, f\,\right > \,+\, \left <\,f,\, e_{\,2}\,\right >\,\left <\,e_{\,2},\, f\,\right >+\cdots\\
&=\,\left <\,f,\, e_{\,1}\,\right >\,\left <\,e_{\,1},\, f\,\right > \,+\, \left <\,f,\, \left <\,f,\, e_{\,1}\,\right >\,e_{\,1}\,\right >\,+\, \left <\,f,\, \left <\,f,\, e_{\,2}\,\right >\,e_{\,2}\,\right >\,+\,\cdots\\
&=\,\left|\,\left <\,f,\, e_{\,1}\,\right >\,\right|^{\,2} \,+\, \left <\,f,\, \sum\limits_{i \,=\, 1}^{\infty}\,\left <\,f,\, e_{\,i}\,\right >\,e_{\,i}\,\right >\\
&\,=\,\left|\,\left <\,f,\, e_{\,1}\,\right >\,\right|^{\,2} \,+\, \left<\,f,\,f\,\right> \,=\, \left|\,\left <\,f,\, e_{\,1}\,\right >\,\right|^{\,2} \,+\, \|\,f\,\|^{\,2}. 
\end{align*} 
Thus, \,$\left(\,F,\, G\,\right)$\, is a continuous biframe for \,$H$\, with bounds \,$1$\, and \,$2$.
\end{example}

Next, we give an example of a continuous biframe for a real inner product space.

\begin{example}
Let
\[
V \,=\, \left\{
\begin{pmatrix}
\;a & 0\\
\;0           & b\\
\end{pmatrix}
:\, a,\, b \,\in\, \mathbb{R}
\right\}.
\]
Define \,$\left<\,\cdot,\, \cdot\,\right> \,:\, V \,\times\, V \,\to\, \mathbb{R}$\, by \,$\left<\,M,\, N\,\right> \,=\, det\left(\,MN^{\,t}\,\right)$, for all \,$M,\, N \,\in V$.\,Then it is easy to verify that \,$\left<\,\cdot,\, \cdot\,\right>$\, is an real inner product on \,$V$.\,Now, we consider a measure space \,$\left(\,\Omega \,=\, [\,0,\,1\,],\, \mu\,\right)$\, where \,$\mu$\,  is the Lebesgue measure.\,Define \,$F \,:\, \Omega \,\to\, V$\, by
\[
F\,(\,w\,) \,=\, 
\begin{pmatrix}
\;\sqrt{3}\,(\,1 \,-\, w\,) & 0\\
\;0           & \sqrt{3}\,(\,1 \,+\, w\,)\\
\end{pmatrix}
,\; w \,\in\, \Omega
\] 
and \,$G \,:\, \Omega \,\to\, V$\, by
\[
G\,(\,w\,) \,=\, 
\begin{pmatrix}
\;\sqrt{2\,w}\, & 0\\
\;0           & \sqrt{2\,w}\,\\
\end{pmatrix}
,\; w \,\in\, \Omega
\]
It is easy to verify that for all \,$M \,\in\, V$, \,$w \,\mapsto\, \left<\,M,\, F\,(\,w\,)\,\right>$\, and \,$w \,\mapsto\, \left<\,M,\, G\,(\,w\,)\,\right>$\, are measurable functions on \,$\Omega$.\,Now, for each \,$
M \,=\, 
\begin{pmatrix}
\;a\, & 0\\
\;0           & b\,\\
\end{pmatrix}
 \,\in\, V
$, we get
\begin{align*}
\left<\,M,\, F\,(\,w\,)\,\right> \,=\, det\left\{\,\begin{pmatrix}
\;a\, & 0\\
\;0           & b\,\\
\end{pmatrix}
\begin{pmatrix}
\;\sqrt{3}\,(\,1 \,-\, w\,) & 0\\
\;0           & \sqrt{3}\,(\,1 \,+\, w\,)\\
\end{pmatrix}
\,\right\} \,=\, 3\,a\,b\,\left(\,1 \,-\, w^{\,2}\,\right)
\end{align*} 
and
\begin{align*}
\left<\,G\,(\,w\,),\, M\,\right> \,=\, det\left\{\,\begin{pmatrix}
\;\sqrt{2\,w}\, & 0\\
\;0           & \sqrt{2\,w}\,\\
\end{pmatrix}
\begin{pmatrix}
\;a\, & 0\\
\;0           & b\,\\
\end{pmatrix}
\,\right\} \,=\, 2\,a\,b\,w.
\end{align*} 
Thus, for each \,$
M \,=\, 
\begin{pmatrix}
\;a\, & 0\\
\;0           & b\,\\
\end{pmatrix}
 \,\in\, V
$, we have
\begin{align*}
&\int\limits_{\,[\,0,\, 1\,]}\,\left<\,M,\, F\,(\,w\,)\,\right>\,\left<\,G\,(\,w\,),\, M\,\right>\,d\mu \,=\, \int\limits_{\,[\,0,\, 1\,]}\,6\,a^{\,2}\,b^{\,2}\,w\,\left(\,1 \,-\, w^{\,2}\,\right)\,d\mu \\
&\hspace{2.5cm}=\,\dfrac{3}{2}\,a^{\,2}\,b^{\,2} \,=\, \dfrac{3}{2}\,\,det \begin{pmatrix}
\;a^{\,2}\, & 0\\
\;0           & b^{\,2}\,\\
\end{pmatrix}
\,=\, \dfrac{3}{2}\, \left\|\,M\,\right\|^{\,2}.
\end{align*}
Therefore, \,$(\,F,\, G\,)$\, is a tight continuous biframe for \,$H$\, with bound \,$3 \,/\, 2$. 
\end{example}

Let \,$(\,F,\, G\,) \,=\, \,\left(\,F \,:\, \Omega \,\to\, H,\, \,G \,:\, \Omega \,\to\, H\,\right)$\, be a continuous biframe for \,$H$\, with respect to \,$\left(\,\Omega,\, \mu\,\right)$.\,Then the mapping \,$\Psi \,:\, H \,\times\, H \,\to\, \mathbb{C}$\, defined by 
\[\Psi\,(\,f,\, g\,) \,=\, \int\limits_{\,\Omega}\,\left<\,f, F\,(\,w\,)\,\right>\left<\,G\,(\,w\,),\, g\,\right>\,d\mu\]
is well-defined, a sesquilinear form (\,i.\,e linear in the first and conjugate-linear in the second variable\,) and is bounded.\,By Cauchy-Schwartz inequality, we get
\begin{align*}
\left|\,\Psi\,(\,f,\, g\,)\,\right| &\,\leq\, \int\limits_{\,\Omega}\,\left|\,\left<\,f, F\,(\,w\,)\,\right>\left<\,G\,(\,w\,),\, g\,\right>\,\right|\,d\mu \\
&\leq\,\left(\,\int\limits_{\,\Omega}\,\left|\,\left<\,f, F\,(\,w\,)\,\right>\,\right|^{\,2}\,d\mu\,\right)^{1 \,/\, 2}\,\left(\,\int\limits_{\,\Omega}\,\left|\,\left<\,G\,(\,w\,),\, g\,\right>\,\right|^{\,2}\,d\mu\,\right)^{1 \,/\, 2}\\
&\leq\,B\,\|\,f\,\|\,\|\,g\,\|
\end{align*}
By Theorem 2.3.6 in \cite{GJ}, there exists a unique operator \,$S_{F,\,G} \,:\, H \,\to\, H$\, such that 
\[\Psi\,(\,f,\, g\,) \,=\, \left<\,S_{F,\,G}\,f,\, g\,\right>\; \;\forall\, f,\, g \,\in\, H\]
and moreover \,$\|\,\Psi\,\| \,=\, \left\|\,S_{F,\,G}\,\right\|$.

Now, we introduce the continuous biframe operator and give some properties.

\begin{definition}
Let \,$(\,F,\, G\,) \,=\, \,\left(\,F \,:\, \Omega \,\to\, H,\, \,G \,:\, \Omega \,\to\, H\,\right)$\, be a continuous biframe for \,$H$\, with respect to \,$\left(\,\Omega,\, \mu\,\right)$.\,Then the operator \,$S_{F,\,G} \,:\, H \,\to\, H$\, defined by
\[S_{F,\,G}\,f \,=\, \int\limits_{\,\Omega}\,\left<\,f, F\,(\,w\,)\,\right>\,G\,(\,w\,)\,d\mu,\]
for all \,$f \,\in\, H$\, is called the frame operator.
\end{definition}

Now, for each \,$f \,\in\, H$, we have 
\begin{align}
\left<\,S_{F,\,G}\,f,\, f\,\right> \,=\, \int\limits_{\,\Omega}\,\left<\,f, F\,(\,w\,)\,\right>\left<\,G\,(\,w\,),\, f\,\right>\,d\mu.\label{3.eqq3.12}
\end{align}
Thus, for each \,$f \,\in\, H$, we get
\[A\,\left\|\,f\,\right\|^{\,2} \leq\, \left<\,S_{F,\,G}\,f,\, f\,\right> \,\leq\,  B\left\|\,f\,\right\|^{\,2}.\]
This implies that \,$A\,I \,\leq\, S_{F,\,G} \,\leq\, B\,I$, where \,$I$\, is the identity operator on \,$H$.\,Hence, \,$S_{F,\,G}$\, is positive and invertible.\,Here, we assume that  \,$S_{F,\,G}$\, is self-adjoint operator.

Thus, every \,$f \,\in\, H$\, has the representations
\begin{align*}
&f \,=\, S_{F,\,G}\, S^{\,-\, 1}_{F,\,G}\,f \,=\, \int\limits_{\,\Omega}\,\left<\,f, S^{\,-\, 1}_{F,\,G}\,F\,(\,w\,)\,\right>\,G\,(\,w\,)\,d\mu\,,\\
&f \,=\, S^{\,-\, 1}_{F,\,G}\,S_{F,\,G}\,f \,=\, \int\limits_{\,\Omega}\,\left<\,f, F\,(\,w\,)\,\right>\,S^{\,-\, 1}_{F,\,G}\,G\,(\,w\,)\,d\mu\,.
\end{align*}

\begin{example}
Let \,$H \,=\, \mathbb{R}^{\,3}$.Suppose
\begin{align*}
&\left\{\,f_{\,i}\,\right\}_{i \,=\, 1}^{\,3} \,=\, \left\{\,(\,2,\, 1,\, 1\,),\, (\,-\, 1,\, 3,\, -\, 1\,),\, (\,-\, 1,\, 1,\, 4\,)\,\right\}\,,\\
&\left\{\,g_{\,i}\,\right\}_{i \,=\, 1}^{\,3} \,=\, \left\{\,(\,1,\, 0,\, 0\,),\, (\,0,\, 1,\, 0\,),\, (\,0,\, 0,\, 1\,)\,\right\}.
\end{align*}
Consider
\[\Omega \,=\, \left\{\,x \,\in\, \mathbb{R}^{\,3} \,:\, \|\,x\,\| \,\leq\, 1\,\right\}.\]
Then \,$\left(\,\Omega,\, \mu\,\right)$\, is a measure space, where \,$\mu$\, is the Lebesgue measure.\,Suppose \\$\{\,B_{\,1},\,B_{\,2},\, B_{\,3}\,\}$\, is a partition of \,$\Omega$\, where \,$\mu\,(\,B_{1}\,) \,\geq\, \mu\,(\,B_{2}\,) \,\geq\, \mu\,(\,B_{3}\,) \,>\, 1$. \\Now, we define 
\[F \,:\, \Omega \,\to\, H\hspace{.5cm}\text{by} \hspace{.3cm} F\,(\,w\,) \,=\, \begin{cases}
\dfrac{f_{1}}{\sqrt{\,\mu\left(\,B_{1}\,\right)}} & \text{if\;\;}\; w \,\in\, B_{1} \\ \dfrac{f_{2}}{\sqrt{\,\mu\left(\,B_{2}\,\right)}} & \text{if\;\;}\; w \,\in\, B_{2}\\ \dfrac{f_{3}}{\sqrt{\,\mu\left(\,B_{3}\,\right)}} & \text{if\;\;}\; w \,\in\, B_{3} \end{cases}\]and
\[G \,:\, \Omega \,\to\, H\hspace{.5cm}\text{by} \hspace{.3cm} G\,(\,w\,) \,=\, \begin{cases}
\dfrac{g_{1}}{\sqrt{\,\mu\left(\,B_{1}\,\right)}} & \text{if\;\;}\; w \,\in\, B_{1} \\ \dfrac{g_{2}}{\sqrt{\,\mu\left(\,B_{2}\,\right)}} & \text{if\;\;}\; w \,\in\, B_{2}\\ \dfrac{g_{3}}{\sqrt{\,\mu\left(\,B_{3}\,\right)}} & \text{if\;\;}\; w \,\in\, B_{3} \end{cases}\]
It is easy to verify that for all \,$f \,\in\, H$, \,$w \,\mapsto\, \left<\,f,\, F\,(\,w\,)\,\right>$\, and \,$w \,\mapsto\, \left<\,f,\, G\,(\,w\,)\,\right>$\, are measurable functions on \,$\Omega$.\,Now, for \,$f \,=\, \left(\,x,\, y,\, z\,\right) \,\in\, H$, we have
\begin{align*}
&\int\limits_{\,\Omega}\,\left<\,f,\, F\,(\,w\,)\,\right>\,\left<\,G\,(\,w\,),\, f\,\right>\,d\mu \\
&= \int\limits_{\,B_{1}}\left<\,f,\, \dfrac{f_{\,1}}{\sqrt{\mu\left(\,B_{1}\,\right)}}\,\right>\,\left<\,\dfrac{g_{\,1}}{\sqrt{\mu\left(\,B_{1}\,\right)}},\, f\,\right>\,d\mu\\
& \,+\, \int\limits_{\,B_{2}}\left<\,f,\, \dfrac{f_{\,2}\,}{\sqrt{\mu\left(\,B_{2}\,\right)}}\,\right>\,\left<\,\dfrac{\, g_{\,2}}{\sqrt{\mu\left(\,B_{2}\,\right)}},\, f\,\right>\,d\mu\\
&+\,\int\limits_{\,B_{3}}\left<\,f,\, \dfrac{f_{\,3}\,}{\sqrt{\mu\left(\,B_{3}\,\right)}}\,\right>\,\left<\, \dfrac{g_{\,3}\,}{\sqrt{\mu\left(\,B_{3}\,\right)}},\, f\,\right>\,d\mu\\
&=\,\left<\,f,\, f_{1}\,\right>\,\left<\,g_{1},\, f\,\right> \,+\, \left<\,f,\, f_{2}\,\right>\,\left<\,g_{2},\, f\,\right> \,+\, \left<\,f,\, f_{3}\,\right>\,\left<\,g_{3},\, f\,\right>\\
&=\,\left<\,\left(\,x,\, y,\, z\,\right),\, (\,2,\, 1,\, 1\,)\,\right>\,\left<\,(\,1,\, 0,\, 0\,),\, \left(\,x,\, y,\, z\,\right)\,\right> \,+\,\\
&\hspace{1cm}+\, \left<\,\left(\,x,\, y,\, z\,\right),\, (\,-\, 1,\, 3,\, \,-\, 1\,)\,\right>\,\left<\,(\,0,\, 1,\, 0\,),\, \left(\,x,\, y,\, z\,\right)\,\right> \\
&\hspace{1cm}+\,\left<\,\left(\,x,\, y,\, z\,\right),\, (\,-\, 1,\, 1,\, 4\,)\,\right>\,\left<\,(\,0,\, 0,\, 1\,),\, \left(\,x,\, y,\, z\,\right)\,\right>\\
&=\,\left(\,2\,x \,+\, y \,+\, z\,\right)\,x \,+\, \left(\,-\, x \,+\, 3\,y \,-\, z\,\right)\,y \,+\, \left(\,-\, x \,+\, y \,+\, 4\,z\,\right)\,z\\
&=\,2\,x^{\,2} \,+\, 3\,y^{\,2} \,+\, 4\,z^{\,2} \,\leq\, 4\,\left(\,x^{\,2} \,+\, y^{\,2} \,+\, z^{\,2}\,\right) \,=\, 4\,\left\|\,\left(\,x,\, y,\, z\,\right)\,\right\|^{2} \,=\, 4\,\left\|\,f\,\right\|^{2}. 
\end{align*} 
Thus, for each \,$f \,\in\, H$, we get
\begin{align*}
2\,\left\|\,f\,\right\|^{2} \,\leq\, \int\limits_{\,\Omega}\,\left<\,f,\, F\,(\,w\,)\,\right>\,\left<\,G\,(\,w\,),\, f\,\right>\,d\mu \,\leq\, 4\,\left\|\,f\,\right\|^{2}.
\end{align*} 
Therefore, \,$(\,F,\, G\,)$\, is a continuous biframe for \,$H$\, with bounds \,$2$\, and \,$4$. 

For \,$\left(\,x,\, y,\, z\,\right) \,\in\, \mathbb{R}^{\,3}$, the continuous biframe operator \,$S_{F,\,G}$\, is given by
\begin{align*}
S_{F,\,G}\left(\,x,\, y,\, z\,\right)& \,=\, \left<\,\left(\,x,\, y,\, z\,\right),\, (\,2,\, 1,\, 1\,)\,\right>\,(\,1,\, 0,\, 0\,) \,+\,\\
&+\, \left<\,\left(\,x,\, y,\, z\,\right),\, (\,-\, 1,\, 3,\, \,-\, 1\,)\,\right>\,(\,0,\, 1,\, 0\,) \\
&+\,\left<\,\left(\,x,\, y,\, z\,\right),\, (\,-\, 1,\, 1,\, 4\,)\,\right>\,(\,0,\, 0,\, 1\,)\\
&=\,\left(\,2\,x \,+\, y \,+\, z,\, \,-\, x \,+\, 3\,y \,-\, z,\, \,-\, x \,+\, y \,+\, 4\,z\,\right).
\end{align*}
The matrix associated with the operator \,$S_{F,\,G}$\, is given by
\[
\left[\,S_{F,\,G}\,\right] \,=\, 
\begin{pmatrix}
\;2 & 1 & 1\\
\;-\, 1   & 3 & \,-\,1\\
\;-\, 1    & 1  & 4
\end{pmatrix}
.\] 
Since \,$det\left(\,\left[\,S_{F,\,G}\,\right]\,\right) \,=\, 33 \,\neq\, 0$, the matrix \,$\left[\,S_{F,\,G}\,\right]$\, is invertible.\,Thus, the operator \,$S_{F,\,G}$\, is well defined and invertible bounded linear operator on \,$\mathbb{R}^{\,3}$.\,It is easy to verify that the operator \,$S_{F,\,G}$\, is positive.       

The inverse of the matrix \,$\left[\,S_{F,\,G}\,\right]$\, is given by 
\[
\left[\,S_{F,\,G}\,\right]^{\,-\, 1} \,=\, \dfrac{1}{33} 
\begin{pmatrix}
\;13 & \,-\,3 & \,-\,4\\
\;5   & 9 &  1\\
\;2    & \,-\,3  & 7
\end{pmatrix}
.\] 
Therefore, for \,$\left(\,x,\, y,\, z\,\right) \,\in\, \mathbb{R}^{\,3}$,  \,$S^{\,-\, 1}_{F,\,G}$\, is given by
\begin{align*}
&S^{\,-\, 1}_{F,\,G}\left(\,x,\, y,\, z\,\right) \,=\, \dfrac{1}{33}\left(\,13\,x \,-\, 3\,y \,-\, 4\,z,\,  5\,x \,+\, 9\,y \,+\, z,\,  2\,x \,-\, 3\,y \,+\, 7\,z\,\right).
\end{align*}
Now, for \,$f \,=\, \left(\,x,\, y,\, z\,\right) \,\in\, H$, we have
\begin{align*}
&\int\limits_{\,\Omega}\,\left<\,f,\, F\,(\,w\,)\,\right>\,S^{\,-\, 1}_{F,\,G}\,G\,(\,w\,)\,d\mu \\
&=\,\left<\,f,\, f_{1}\,\right>\,S^{\,-\, 1}_{F,\,G}\,g_{1} \,+\, \left<\,f,\, f_{2}\,\right>\,S^{\,-\, 1}_{F,\,G}\,g_{2} \,+\, \left<\,f,\, f_{3}\,\right>\,S^{\,-\, 1}_{F,\,G}\,g_{3}\\
&=\,\left(\,2\,x \,+\, y \,+\, z\,\right)\,\dfrac{1}{33}\,\left(\,13,\, 5,\, 2\,\right) \,+\, \left(\,-\, x \,+\, 3\,y \,-\, z\,\right)\,\dfrac{1}{33}\,\left(\,-\, 3,\, 9,\, \,-\, 3\,\right) \,+\\
&\hspace{1cm}\,+\, \left(\,-\, x \,+\, y \,+\, 4\,z\,\right)\,\dfrac{1}{33}\,\left(\,-\,4,\, 1,\, 7\,\right)\\
&=\,\dfrac{1}{33}\,\left(\,33\,x,\, 33\,y,\, 33\,z\,\right) \,=\, \left(\,x,\, y,\, z\,\right) \,=\, f
\end{align*} 
Similarly, it can be verified that
\begin{align*}
&\int\limits_{\,\Omega}\,\left<\,f,\, S^{\,-\, 1}_{F,\,G}\,F\,(\,w\,)\,\right>\,G\,(\,w\,)\,d\mu \,=\, f\,,\, f \,\in\, H. 
\end{align*}
Thus,the representation theorem is verified in this example.
\end{example}

\begin{theorem}
The pair \,$(\,F,\, G\,)$\, is a continuous biframe for \,$H$\, with respect to \,$\left(\,\Omega,\, \mu\,\right)$\, if and only if  \,$(\,G,\, F\,)$\, is a continuous biframe for \,$H$\, with respect to \,$\left(\,\Omega,\, \mu\,\right)$.
\end{theorem}

\begin{proof}
Let \,$(\,F,\, G\,)$\, is a continuous biframe for \,$H$\, with bounds \,$A$\, and \,$B$.\,Then for each \,$f \,\in\, H$, we have 
\[A\,\left\|\,f\,\right\|^{\,2} \leq \int\limits_{\Omega}\,\left<\,f,\, F\,(\,w\,)\,\right>\,\left<\,G\,(\,w\,),\, f\,\right>\,d\mu \leq\, B\left\|\,f\,\right\|^{\,2}.\] 
Since \,$S_{F,\,G}$\, is self adjoint, using (\ref{3.eqq3.12}), we can write
\begin{align*}
\int\limits_{\Omega}\,\left<\,f,\, F\,(\,w\,)\,\right>\,\left<\,G\,(\,w\,),\, f\,\right>\,d\mu & \,=\, \left<\,S_{F,\,G}\,f,\, f\,\right> \,=\, \overline{\left<\,S_{F,\,G}\,f,\, f\,\right>}\\  
&\,=\, \overline{\int\limits_{\Omega}\,\left<\,f,\, F\,(\,w\,)\,\right>\,\left<\,G\,(\,w\,),\, f\,\right>\,d\mu}\\
&=\, \int\limits_{\Omega}\,\left<\,f,\, G\,(\,w\,)\,\right>\,\left<\,F\,(\,w\,),\, f\,\right>\,d\mu.
\end{align*}
Thus, for each \,$f \,\in\, H$, we have 
\[A\,\left\|\,f\,\right\|^{\,2} \leq \int\limits_{\Omega}\,\left<\,f,\, G\,(\,w\,)\,\right>\,\left<\,F\,(\,w\,),\, f\,\right>\,d\mu \leq\, B\,\left\|\,f\,\right\|^{\,2}.\]
Therefore, \,$(\,G,\, F\,)$\, is a continuous biframe for \,$H$.

Similarly, we can prove the converse part of this Theorem.  
\end{proof}

In the next Theorem, we establish a characterization of a continuous biframe using its biframe operator.

\begin{theorem}
Let \,$(\,F,\, G\,)$\, is a continuous biframe Bessel mapping for \,$H$\, with respect to \,$\left(\,\Omega,\, \mu\,\right)$.\,Then \,$(\,F,\,G\,)$\, is a continuous biframe for \,$H$\, if and only if there exists \,$\alpha \,>\, 0$\, such that \,$S_{F,\,G} \,\geq\, \alpha\,I$, where \,$S_{F,\,G}$\, is the continuous biframe operator for \,$(\,F,\, G\,)$.
\end{theorem}

\begin{proof}
Let \,$(\,F,\, G\,)$\, is a continuous biframe for \,$H$\, with bounds \,$A$\, and \,$B$.\,Then using (\ref{3.eqq3.11}) and (\ref{3.eqq3.12}), for each \,$f \,\in\, H$, we get 
\[A\,\left\|\,f\,\right\|^{\,2} \leq\, \left<\,S_{F,\,G}\,f,\, f\,\right> \,\leq\,  B\left\|\,f\,\right\|^{\,2}.\]
Thus
\[A\,\left<\,f,\, f\,\right> \leq\, \left<\,S_{F,\,G}\,f,\, f\,\right> \,\Rightarrow\, S_{F,\,G} \,\geq\, \alpha\,I,\]
where \,$\alpha \,=\, A$.

Conversely, suppose that \,$S_{F,\,G} \,\geq\, \alpha\,I$.\,Thus, for each \,$f \,\in\, H$, we have
\[\alpha^{\,2}\,\|\,f\,\|^{\,2} \,\leq\,  \left<\,S_{F,\,G}\,f,\, f\,\right> \,=\, \int\limits_{\Omega}\,\left<\,f,\, F\,(\,w\,)\,\right>\,\left<\,G\,(\,w\,),\, f\,\right>\,d\mu.\] 
Hence, \,$(\,F,\, G\,)$\, is a continuous biframe for \,$H$\, with lower biframe bound \,$\alpha^{\,2}$.\,This completes the proof.
\end{proof}

Next, also we give a characterization of a continuous biframe with the help of a invertible operator on \,$H$.

\begin{theorem}\label{3.thm3.39}
Let \,$T$\, be an invertible bounded linear operator on \,$H$.\,Then \,$(\,F,\, G\,)$\, is a continuous biframe for \,$H$\, with respect to \,$\left(\,\Omega,\, \mu\,\right)$\, if and only if \,$(\,TF,\, TG\,)$\, is a continuous biframe for \,$H$\, with respect to \,$\left(\,\Omega,\, \mu\,\right)$.
\end{theorem}

\begin{proof}
For each \,$f \,\in\, H$, \,$w \,\mapsto\, \left<\,f,\, TF\,(\,w\,)\,\right> \,=\, \left<\,T^{\,\ast}\,f,\, F\,(\,w\,)\,\right>$\, and \,$w \,\mapsto\, \left<\,f,\, TG\,(\,w\,)\,\right> \,=\, \left<\,T^{\,\ast}\,f,\, G\,(\,w\,)\,\right>$\, are measurable functions on \,$\Omega$.\,Let \,$(\,F,\, G\,)$\, is a continuous biframe for \,$H$\, with bounds \,$A$\, and \,$B$.\,Since \,$T$\, is invertible, for \,$f \,\in\, H$, we have
\begin{align*}
&\left\|\,f\,\right\|^{\,2} \,=\, \left\|\,\left (\,T^{\,-\, 1}\,\right )^{\,\ast}\,T^{\,\ast}\,f\,\right\|^{\,2} \,\leq\, \left\|\,T^{\,-\, 1}\,\right\|^{\,2} \,\left\|\,T^{\,\ast}\,f\,\right \|^{\,2}.
\end{align*}
Now, for each \,$f \,\in\, H$, we have 
\begin{align*}
\int\limits_{\Omega}\,\left<\,f,\, TF\,(\,w\,)\,\right>\,\left<\,TG\,(\,w\,),\, f\,\right>\,d\mu &\,=\, \int\limits_{\Omega}\,\left<\,T^{\,\ast}\,f,\, F\,(\,w\,)\,\right>\,\left<\,G\,(\,w\,),\, T^{\,\ast}\,f\,\right>\,d\mu\\
&\leq\,B\,\left\|\,T^{\,\ast}\,f\,\right\|^{\,2}\,\leq\,B\,\|\,T\,\|^{\,2}\,\left\|\,f\,\right\|^{\,2}.
\end{align*}
On the other hand, for each \,$f \,\in\, H$, we have 
\begin{align*}
\int\limits_{\Omega}\,\left<\,f,\, TF\,(\,w\,)\,\right>\,\left<\,TG\,(\,w\,),\, f\,\right>\,d\mu &\,=\, \int\limits_{\Omega}\,\left<\,T^{\,\ast}\,f,\, F\,(\,w\,)\,\right>\,\left<\,G\,(\,w\,),\, T^{\,\ast}\,f\,\right>\,d\mu\\
&\geq\,A\,\left\|\,T^{\,\ast}\,f\,\right\|^{\,2}\,\geq\,A\,\left\|\,T^{\,-\, 1}\,\right\|^{\,-\, 2}\,\left\|\,f\,\right\|^{\,2}.
\end{align*}
Hence, \,$(\,TF,\, TG\,)$\, is a continuous biframe for \,$H$\, with bounds \,$A\,\left\|\,T^{\,-\, 1}\,\right\|^{\,-\, 2}$\, and \,$B\,\|\,T\,\|^{\,2}$.

Conversely, suppose that \,$(\,TF,\, TG\,)$\, is a continuous biframe for \,$H$\, with bounds \,$A$\, and \,$B$.\,Now, for each \,$f \,\in\, H$, we have
\begin{align*}
&\dfrac{A}{\|\,T\,\|^{\,2}}\,\left\|\,f\,\right\|^{\,2} \,=\, \dfrac{A}{\|\,T\,\|^{\,2}}\,\left\|\,\left(\,T^{\,-\, 1}\,T\,\right)^{\,\ast}\,f\,\right\|^{\,2} \,\leq\, A\,\left\|\,\left(\,T^{\,-\, 1}\,\right)^{\,\ast}\,f\,\right\|^{\,2} \\
&\leq\, \int\limits_{\Omega}\,\left<\,\left(\,T^{\,-\, 1}\,\right)^{\,\ast}\,f,\, TF\,(\,w\,)\,\right>\,\left<\,TG\,(\,w\,),\, \left(\,T^{\,-\, 1}\,\right)^{\,\ast}\,f\,\right>\,d\mu \\
&\,=\, \int\limits_{\Omega}\,\left<\,T^{\,\ast}\left(\,T^{\,-\, 1}\,\right)^{\,\ast}\,f,\, F\,(\,w\,)\,\right>\,\left<\,G\,(\,w\,),\, T^{\,\ast}\left(\,T^{\,-\, 1}\,\right)^{\,\ast}\,f\,\right>\,d\mu\\
&\,=\, \int\limits_{\Omega}\,\left<\,f,\, F\,(\,w\,)\,\right>\,\left<\,G\,(\,w\,),\, f\,\right>\,d\mu.   
\end{align*}
On the other hand, for each \,$f \,\in\, H$, we have
\begin{align*}
&\int\limits_{\Omega}\,\left<\,f,\, F\,(\,w\,)\,\right>\,\left<\,G\,(\,w\,),\, f\,\right>\,d\mu\\
&=\,\int\limits_{\Omega}\,\left<\,T^{\,\ast}\left(\,T^{\,-\, 1}\,\right)^{\,\ast}\,f,\, F\,(\,w\,)\,\right>\,\left<\,G\,(\,w\,),\, T^{\,\ast}\left(\,T^{\,-\, 1}\,\right)^{\,\ast}\,f\,\right>\,d\mu\\
&=\,\int\limits_{\Omega}\,\left<\,\left(\,T^{\,-\, 1}\,\right)^{\,\ast}\,f,\, TF\,(\,w\,)\,\right>\,\left<\,TG\,(\,w\,),\, \left(\,T^{\,-\, 1}\,\right)^{\,\ast}\,f\,\right>\,d\mu\\
&\leq\,B\,\left\|\,\left(\,T^{\,-\, 1}\,\right)^{\,\ast}\,f\,\right\|^{\,2} \,\leq\, B\,\left\|\,T^{\,-\, 1}\,\right\|^{\,2}\,\left\|\,f\,\right\|^{\,2}.
\end{align*}
Thus, \,$(\,F,\, G\,)$\, is a continuous biframe for \,$H$\, with bounds \,$\dfrac{A}{\|\,T\,\|^{\,2}}$\, and \,$B\,\left\|\,T^{\,-\, 1}\,\right\|^{\,2}$.\,This completes the proof. 
\end{proof}

Now, we would complete this section with discussion of continuous biframe Bessel multiplier in \,$H$.

\begin{definition}
Let \,$(\,F,\, F\,)$\, and \,$(\,G,\, G\,)$\, be continuous biframe Bessel mappings for \,$H$\, with respect to \,$\left(\,\Omega,\, \mu\,\right)$\, and let \,$m \,:\, \Omega \,\to\, \mathbb{C}$\, be a measurable function.\,Then the operator \,$M_{m,\, F,\, G} \,:\, H \,\to\, H$\, defined by
\begin{align*}
&\left<\,M_{m,\, F,\, G}\,f,\, g \,\right> \,=\, \int\limits_{\,\Omega}\,m\,(\,w\,)\,\left<\,f,\, F\,(\,w\,) \,\right>\,\left<\,G\,(\,w\,),\, g\,\right>\,d\mu,
\end{align*} 
for all \,$f,\, g \,\in\, H$, is called continuous biframe Bessel multiplier of \,$F$\, and \,$G$\, with respect to \,$m$. 
\end{definition}

\begin{theorem}\label{thm4.22}
The continuous biframe Bessel multiplier of \,$F$\, and \,$G$\, with respect to \,$m$\, is well defined and bounded. 
\end{theorem}

\begin{proof}
Let \,$(\,F,\, F\,)$\, and \,$(\,G,\, G\,)$\, be continuous biframe Bessel mappings for \,$H$\, with bounds \,$B_{1}$\, and \,$B_{2}$.\,Then for any \,$f,\, g \,\in\, H$, we have
\begin{align*}
&\left|\,\left<\,M_{m,\, F,\, G}\,f,\, g\,\right>\,\right| =\,\left|\int\limits_{\,\Omega}\,m\,(\,w\,)\left<\,f,\, F\,(\,w\,)\,\right>\left<\,G\,(\,w\,),\, g\,\right>d\mu\,\right|\\
&\leq\,\|\,m\,\|_{\,\infty}\left(\,\int\limits_{\,\Omega}\,\left|\,\left<\,f,\, F\,(\,w\,)\,\right>\,\right|^{\,2}d\mu\,\right)^{1 \,/\, 2}\,\left(\,\int\limits_{\,\Omega}\,\left|\,\left<\,g,\, G\,(\,w\,)\,\right>\,\right|^{\,2}d\mu\,\right)^{1 \,/\, 2}\\
&\leq\,\|\,m\,\|_{\,\infty}\,\sqrt{B_{\,1}\,B_{\,2}}\,\left\|\,f\,\right\|\,\left\|\,g\,\right\|.
\end{align*}
This shows that \,$\left\|\,M_{m,\, F,\, G}\,\right\| \,\leq\, \|\,m\,\|_{\,\infty}\,\sqrt{B_{\,1}\,B_{\,2}}$\, and so \,$M_{m,\, F,\, G}$\, is well-defined and bounded.\,This completes the proof. 
\end{proof}

Following the proof of the Theorem \ref{thm4.22}, for each \,$f \,\in\, H$, we have
\begin{align}
&\left\|\,M_{m,\, F,\, G}\,f\,\right\| \,=\, \sup\limits_{\left\|\,g\,\right\| \,=\, 1}\,\left|\,\left<\,M_{m,\, F,\, G}\,f,\, g\,\right>\,\right|\nonumber\\
&\leq\,\|\,m\,\|_{\,\infty}\,\sqrt{\,B_{2}}\,\left(\,\int\limits_{\,\Omega}\,\left|\,\left<\,f,\, F\,(\,w\,)\,\right>\,\right|^{\,2}d\mu\,\right)^{1 \,/\, 2}\label{en4.22}
\end{align}
and similarly it can be shown that
\begin{align}
&\left\|\,M_{m,\, F,\, G}^{\,\ast}\,g\,\right\|\nonumber\\
& \,\leq\,\|\,m\,\|_{\,\infty}\,\sqrt{\,B_{1}}\,\left(\,\int\limits_{\,\Omega}\,\left|\,\left<\,G\,(\,w\,),\, g\,\right>\,\right|^{\,2}d\mu\,\right)^{1 \,/\, 2}.\label{en4.23} 
\end{align}

\begin{theorem}
Let \,$M_{m,\, F,\, G}$\, be the continuous biframe Bessel multiplier of \,$F$\, and \,$G$\, with respect to \,$m$.\,Then \,$M^{\,\ast}_{m,\, F,\, G} \,=\, M_{\overline{m},\, F,\, G}$.
\end{theorem}

\begin{proof}
For \,$f,\, g \,\in\, H$, we have
\begin{align*}
&\left<\,f,\, M^{\,\ast}_{m,\, F,\, G}\,g\,\right> \,=\, \left<\,M_{m,\, F,\, G}\,f,\, g \,\right>\\
& \,=\, \int\limits_{\,\Omega}\,m\,(\,w\,)\,\left<\,f,\, F\,(\,w\,) \,\right>\,\left<\,G\,(\,w\,),\, g\,\right>\,d\mu\\
&=\,\int\limits_{\,\Omega}\,\left<\,f,\, \overline{m}\,(\,w\,)\left<\,g,\, G\,(\,w\,)\,\right>\,F\,(\,w\,)\,\right>\,\,d\mu\\
&=\, \left<\,f,\, M_{\overline{m},\, F,\,G}\,g \,\right>.
\end{align*}
This completes the proof.
\end{proof}

\begin{theorem}\label{2.thm2.22}
Let \,$M_{m,\, F,\, G}$\, be the continuous biframe Bessel multiplier of \,$F$\, and \,$G$\, with respect to \,$m$.\,Then \,$\left(\,F,\, F\,\right)$\, is a continuous biframe for \,$H$\, provided for each \,$f \,\in\, H$, there exists \,$D \,>\, 0$\, such that
\[\left\|\,M_{m,\, F,\, G}\,f\,\right\| \,\geq\, D\,\left\|\,f\,\right\|.\]
\end{theorem}

\begin{proof}
For each  \,$f \,\in\, H$, using (\ref{en4.22}), we get
\begin{align*}
&D^{\,2}\,\left\|\,f\,\right\|^{\,2}\, \leq\, \left\|\,M_{m,\, F,\, G}\,f\,\right\|^{\,2}\\ 
&\Rightarrow\, D^{\,2}\,\left\|\,f\,\right\|^{\,2}\, \leq\, \|\,m\,\|^{\,2}_{\,\infty}\,B_{\,2}\,\int\limits_{\,\Omega}\,\left|\,\left<\,f,\, F\,(\,w\,)\,\right>\,\right|^{\,2}d\mu\\
&\Rightarrow\,\dfrac{D^{\,2}}{\|\,m\,\|^{\,2}_{\,\infty}\,B_{\,2}}\,\left\|\,f\,\right\|^{\,2}\, \leq\, \int\limits_{\,\Omega}\,\left|\,\left<\,f,\, F\,(\,w\,)\,\right>\,\right|^{\,2}d\mu. \\
&\Rightarrow\,\dfrac{D^{\,2}}{\|\,m\,\|^{\,2}_{\,\infty}\,B_{\,2}}\,\left\|\,f\,\right\|^{\,2}\, \leq\, \int\limits_{\,\Omega}\,\left<\,f,\, F\,(\,w\,)\,\right>\,\left<\,F\,(\,w\,),\, f\,\right>\,d\mu.
\end{align*} 
Thus,  \,$\left(\,F,\, F\,\right)$\, is a continuous biframe for \,$H$\, with bounds \,$\dfrac{D^{\,2}}{\|\,m\,\|^{\,2}_{\,\infty}\,B_{\,2}}$\, and \,$B_{1}$.\,This completes the proof.  
\end{proof}

\begin{theorem}\label{2.th2.4.16}
Let \,$M_{m,\, F,\, G}$\, be the continuous biframe Bessel multiplier of \,$F$\, and \,$G$\, with respect to \,$m$.\,Suppose \,$\lambda_{\,1} \,<\, 1,\, \,\lambda_{\,2} \,>\, -\, 1$\, such that for each \,$f \,\in\, H$, we have
\[\left\|\,f \,-\, M_{m,\, F,\, G}\,f\,\right\| \,\leq\, \lambda_{1}\,\left\|\,f\,\right\| \,+\, \lambda_{2}\,\left\|\,M_{m,\, F,\, G}\,f\,\right\|.\]
Then \,$\left(\,F,\, F\,\right)$\, is a continuous biframe for \,$H$.
\end{theorem}

\begin{proof}
For each \,$f \,\in\, H$, we have
\begin{align*}
&\left\|\,f\,\right\| \,-\, \left\|\,M_{m,\, F,\, G}\,f\,\right\|\leq\,\left\|\,f \,-\, M_{m,\, F,\, G}\,f\,\right\|\\
& \,\leq\, \lambda_{1}\,\left\|\,f\,\right\| \,+\, \lambda_{2}\,\left\|\,M_{m,\, F,\, G}\,f\,\right\|\\
&\Rightarrow\,\left(\,1 \,-\, \lambda_{1}\,\right)\,\left\|\,f\,\right\| \,\leq\, \left(\,1 \,+\, \lambda_{2}\,\right)\,\left\|\,M_{m,\, F,\, G}\,f\,\right\|. 
\end{align*}
Now, using (\ref{en4.22}), we get
\begin{align*}
&\left(\,\dfrac{1 \,-\, \lambda_{1}}{1 \,+\, \lambda_{2}}\,\right)\,\left\|\,f\,\right\|\\
&\,\leq\, \|\,m\,\|_{\,\infty}\,\sqrt{\,B_{2}}\,\left(\,\int\limits_{\,\Omega}\,\left|\,\left<\,f,\, F\,(\,w\,)\,\right>\,\right|^{\,2}d\mu\,\right)^{1 \,/\, 2}.
\end{align*}
\begin{align}
&\Rightarrow\,\dfrac{\left(\,1 \,-\, \lambda_{1}\,\right)^{\,2}}{\|\,m\,\|^{\,2}_{\,\infty}\,B_{2}\,\left(\,1 \,+\, \lambda_{2}\,\right)^{\,2}}\,\left\|\,f\,\right\|^{\,2}\nonumber\\
&\hspace{2cm}\leq\,\int\limits_{\,\Omega}\,\left<\,f,\, F\,(\,w\,)\,\right>\,\left<\,F\,(\,w\,),\, f\,\right>\,d\mu.\label{en4.24} 
\end{align}
Thus, \,$\left(\,F,\, F\,\right)$\, is a continuous biframe for \,$H$\, with bounds \,$\dfrac{\left(\,1 \,-\, \lambda_{1}\,\right)^{\,2}}{\|\,m\,\|^{\,2}_{\,\infty}\,B_{2}\,\left(\,1 \,+\, \lambda_{2}\,\right)^{\,2}}$\, and \,$B_{1}$.\,This completes the proof.  
\end{proof}

\begin{theorem}
Let \,$M_{m,\, F,\, G}$\, be the continuous biframe Bessel multiplier of \,$F$\, and \,$G$\, with respect to \,$m$.\,Suppose \,$\lambda \,\in\, [\,0,\,1\,)$\, such that for each \,$f \,\in\, H$, we have
\[\left\|\,f \,-\, M_{m,\, F,\, G}\,f\,\right\| \,\leq\, \lambda\,\left\|\,f\,\right\|.\]
Then \,$\left(\,F,\, G\,\right)$\, and \,$\left(\,G,\, G\,\right)$\, are continuous biframes for \,$H$.
\end{theorem}

\begin{proof}
Putting \,$\lambda_{1} \,=\, \lambda$\, and \,$\lambda_{2} \,=\, 0$\, in (\ref{en4.24}), we get
\[\dfrac{\left(\,1 \,-\, \lambda\,\right)^{\,2}}{\|\,m\,\|^{\,2}_{\,\infty}\,B_{2}\,}\,\left\|\,f\,\right\|^{\,2}\leq\,\int\limits_{\,\Omega}\,\left<\,f,\, F\,(\,w\,)\,\right>\,\left<\,F\,(\,w\,),\, f\,\right>\,d\mu.\]
Thus, \,$\left(\,F,\, F\,\right)$\, is a continuous biframe for \,$H$.\,On the other hand, for each \,$f \,\in\, H$, we have
\begin{align*}
&\left\|\,f \,-\, M_{m,\, F,\, G}^{\,\ast}\,f\,\right\| \,=\, \left\|\,\left(\,I \,-\, M_{m,\, F,\, G}\,\right)^{\,\ast}\,f\,\right\|\\
&\leq\,\left\|\,I \,-\, M_{m,\, F,\, G}\,\right\|\,\left\|\,f\,\right\| \,\leq\, \lambda\,\left\|\,f\,\right\|\\
&\Rightarrow\,\left(\,1 \,-\, \lambda\,\right)\,\left\|\,f\,\right\| \,\leq\, \left\|\,M_{m,\, F,\, G}^{\,\ast}\,f\,\right\|.  
\end{align*}
Now, using (\ref{en4.23}), we get
\[\dfrac{\left(\,1 \,-\, \lambda\,\right)^{\,2}}{\|\,m\,\|^{\,2}_{\,\infty}\,B_{1}\,}\,\left\|\,f\,\right\|^{\,2}\leq\,\int\limits_{\,\Omega}\,\left|\,\left<\,G\,(\,w\,),\, f\,\right>\,\right|^{\,2}d\mu \,=\, \int\limits_{\,\Omega}\,\left<\,f,\, G\,(\,w\,)\,\right>\,\left<\,G\,(\,w\,),\, f\,\right>\,d\mu.\] 
Thus \,$\left(\,G,\, G\,\right)$\, is a continuous biframe for \,$H$\, with bounds \,$\dfrac{\left(\,1 \,-\, \lambda\,\right)^{\,2}}{\|\,m\,\|^{\,2}_{\,\infty}\,B_{1}\,}$\, and \,$B_{2}$.\,This completes the proof.  
\end{proof}

\begin{definition}
Let \,$\left(\,F,\, G\,\right)$\, be a continuous biframe for \,$H$.\,If   
\begin{align*}
&\left<\,f,\,g\,\right> \,=\, \int\limits_{\,\Omega}\,\left<\,f,\, F\,(\,w\,)\,\right>\,\left<\,G\,(\,w\,),\, g\,\right>\,d\mu\,,
\end{align*} 
holds for all \,$f,\, g \,\in\, H$, then \,$\left(\,F,\, G\,\right)$\, is called dual continuous biframe for \,$H$.
\end{definition}  

\begin{theorem}
Let \,$\left(\,F,\, G\,\right)$\, be a continuous biframe for \,$H$\, with continuous biframe operator \,$S_{F,\,G}$.\,Then \,$\left(\,S^{\,-\, 1}_{F,\,G}\,F,\, G\,\right)$\, and \,$\left(\,F,\, S^{\,-\, 1}_{F,\,G}\,G\,\right)$\, are dual continuous biframes for \,$H$.  
\end{theorem}

\begin{proof}
For each \,$f,\, g \,\in\, H$, we have 
\begin{align*}
&\left<\,f,\, g\,\right> \,=\, \int\limits_{\,\Omega}\,\left<\,f, S^{\,-\, 1}_{F,\,G}\,F\,(\,w\,)\,\right>\,\left<\,G\,(\,w\,),\, g\,\right>\,d\mu\,,\\
&\left<\,f,\, g\,\right> \,=\, \int\limits_{\,\Omega}\,\left<\,f, F\,(\,w\,)\,\right>\,\left<\,S^{\,-\, 1}_{F,\,G}\,G\,(\,w\,),\, g\,\right>\,d\mu\,.
\end{align*}
This verifies that \,$\left(\,S^{\,-\, 1}_{F,\,G}\,F,\, G\,\right)$\, and \,$\left(\,F,\, S^{\,-\, 1}_{F,\,G}\,G\,\right)$\, are dual continuous biframes for \,$H$.   
\end{proof}

In the following Theorem, we will find a dual continuous biframe for \,$H$\, with respect to the multiplier operator.

\begin{theorem}
Let \,$M_{m,\, F,\, G}$\, be invertible and \,$\left(\,F,\, G\,\right)$\, be a continuous biframe for \,$H$.\,Then \,$\left(\,\left(\,M^{\,-\, 1}_{m,\, F,\, G}\,\overline{m}\,F\,\right)^{\,\ast},\, G\,\right)$\, is a dual continuous biframe for \,$H$.   
\end{theorem}

\begin{proof}
From the definition of \,$M_{m,\, F,\, G}$, we can write
\begin{align*}
&\left<\,M_{m,\, F,\, G}\,f,\, g\,\right>\\
& \,=\, \int\limits_{\,\Omega}\,m\,(\,w\,)\,\left<\,f,\,F\,(\,w\,) \,\right>\,\left<\,G\,(\,w\,),\, g\,\right>\,d\mu.
\end{align*}
Now, by replacing \,$f$\, with \,$M^{\,-\, 1}_{m,\, F,\, G}\,f$, we get
\begin{align*}
&\left<\,f,\, g\,\right>\\
& \,=\, \int\limits_{\,\Omega}\,m\,(\,w\,)\,\left<\,M^{\,-\, 1}_{m,\, F,\, G}\,f,\, F\,(\,w\,)\,\right>\,\left<\,G\,(\,w\,),\, g\,\right>\,d\mu\\
&=\,\int\limits_{\,\Omega}\,\left<\,f,\, \left(\,M^{\,-\, 1}_{m,\, F,\, G}\,\right)^{\,\ast}\,\overline{m}\,(\,w\,)\,F\,(\,w\,)\,\right>\,\left<\,G\,(\,w\,),\, g\,\right>\,d\mu.  
\end{align*}
Thus, \,$\left(\,\left(\,M^{\,-\, 1}_{m,\, F,\, G}\,\overline{m}\,F\,\right)^{\,\ast},\, G\,\right)$\, is a dual continuous biframe for \,$H$.\,This completes the proof.    
\end{proof}

\section{Continuous biframe in $H_{1} \,\otimes\, H_{2}$}

In this section, we introduce the concept of continuous biframe in tensor product of Hilbert spaces \,$H_{1} \,\otimes\, H_{2}$\, and give a characterization.

\begin{definition}
Let \,$(\,X,\, \mu\,) \,=\, \left(\,X_{1} \,\times\, X_{2},\, \mu_{\,1} \,\otimes\, \mu_{\,2}\,\right)$\, be the product of measure spaces with \,$\sigma$-finite positive measures \,$\mu_{\,1},\, \mu_{\,2}$\, on \,$X_{1}$, \,$X_{2}$, respectively.\,A pair \,$(\,\mathbf{F},\, \mathbf{G}\,) \,=\, \,\left(\,\mathbb{F} \,:\, X \,\to\, H_{1} \,\otimes\, H_{2},\, \,\mathbb{G} \,:\, X \,\to\, H_{1} \,\otimes\, H_{2}\,\right)$\, is called a continuous biframe for \,$H_{1} \,\otimes\, H_{2}$\, with respect to \,$(\,X,\, \mu\,)$\, if
\begin{itemize}
\item[$(i)$]$ \mathbb{F},\, \mathbb{G}$\, is weakly-measurable, i.\,e., for all \,$f \,\otimes\, g \,\in\, H_{1} \,\otimes\, H_{2}$, \,$x \,=\, \left(\,x_{\,1},\, x_{\,2}\,\right) \,\mapsto\, \left<\,f \,\otimes\, g,\, \mathbb{F}\,(\,x\,)\,\right>$\, and \,$\left(\,x_{\,1},\, x_{\,2}\,\right) \,\mapsto\, \left<\,f \,\otimes\, g,\, \mathbb{G}\,(\,x\,)\,\right>$\, are measurable functions on \,$X$,
\item[$(ii)$]there exist constants \,$A,\,B \,>\, 0$\, such that
\begin{align}
&A\,\left\|\,f \,\otimes\, g\,\right\|^{\,2}\nonumber\\
& \,\leq\, \int\limits_{\,X}\,\left<\,f \,\otimes\, g,\, \mathbb{F}\,(\,x\,)\,\right>\,\left<\,\mathbb{G}\,(\,x\,),\, f \,\otimes\, g\,\right>\,d\mu\nonumber\\
&\leq\,B\,\left\|\,f \,\otimes\, g\,\right\|^{\,2},\label{4.eqq4.11} 
\end{align}
for all \,$f \,\otimes\, g \,\in\, H_{1} \,\otimes\, H_{2}$.\,The constants \,$A$\, and \,$B$\, are called continuous biframe bounds.\,If \,$A \,=\, B$, then the pair $(\,\mathbf{F},\, \mathbf{G}\,)$\, is called a tight continuous biframe for \,$H_{1} \,\otimes\, H_{2}$.\,If \,$(\,\mathbf{F},\, \mathbf{G}\,)$\, satisfies only the right inequality of (\ref{4.eqq4.11}), then it is called continuous biframe Bessel mapping in \,$H_{1} \,\otimes\, H_{2}$\, with Bessel bound \,$B$.  
\end{itemize}   
\end{definition}

\begin{theorem}\label{th4.11}
The pair of mappings \,$\left(\,\mathbf{F},\, \mathbf{G}\,\right) \,=\,  \left(\,F_{1} \,\otimes\, F_{2},\, G_{1} \,\otimes\, G_{2}\,\right) \,=\, \mathbb{F},\, \mathbb{G} \,:\, X \,\to\, H_{1} \,\otimes\, H_{2}$\, is a continuous biframe for \,$H_{1} \,\otimes\, H_{2}$\, with respect to \,$(\,X,\, \mu\,)$\, if and only if \,$F_{1},\,G_{1} \,:\, X_{1} \,\to\, H_{1}$\, is a continuous biframe for \,$H_{1}$\, with respect to \,$\left(\,X_{1},\, \mu_{\,1}\,\right)$\, and \,$F_{2},\, G_{2} \,:\, X_{2} \,\to\, H_{2}$\, is a continuous biframe for \,$H_{2}$\, with respect to \,$\left(\,X_{2},\, \mu_{\,2}\,\right)$.  
\end{theorem}
 
\begin{proof}
Suppose that \,$\left(\,\mathbf{F},\, \mathbf{G}\,\right) \,=\,  \left(\,F_{1} \,\otimes\, F_{2},\, G_{1} \,\otimes\, G_{2}\,\right)$\, is a continuous biframe for \,$H_{1} \,\otimes\, H_{2}$\, with respect to \,$(\,X,\, \mu\,)$\, having bounds \,$A$\, and \,$B$.\,Let \,$f \,\in\, H_{1} \,-\, \{\,\theta\,\}$\, and fix \,$g \,\in\, H_{2} \,-\, \{\,\theta\,\}$.\,Then \,$f \,\otimes\, g \,\in\, H_{1} \,\otimes\, H_{2} \,-\, \{\,\theta\,\otimes\,\theta\,\}$\, and by Fubini's theorem we have
\begin{align*}
&\int\limits_{\,X}\,\left<\,f \,\otimes\, g,\, F_{1}\,(\,x_{\,1}\,) \,\otimes\, F_{2}\,(\,x_{\,2}\,)\,\right>\,\left<\, G_{1}\,(\,x_{\,1}\,) \,\otimes\, G_{2}\,(\,x_{\,2}\,),\, f \,\otimes\, g\,\right>\,d\mu\\
&=\,\int\limits_{\,X_{1}}\,\left<\,f,\, F_{1}\,(\,x_{\,1}\,)\,\right>_{1}\,\left<\,G_{1}\,(\,x_{\,1}\,),\, f\,\right>_{1}\,d\mu_{\,1}\,\int\limits_{\,X_{2}}\,\left<\,g,\, F_{2}\,(\,x_{\,2}\,)\,\right>_{2}\,\left<\,G_{2}\,(\,x_{\,2}\,),\, g\,\right>_{2}\,d\mu_{\,2}.
\end{align*}
Therefore, for each \,$f \,\otimes\, g \,\in\, H_{1} \,\otimes\, H_{2}$, the inequality (\ref{4.eqq4.11}) can be written as
\begin{align*}
&A\,\left\|\,f\,\right\|^{\,2}_{1}\,\left\|\,g\,\right\|^{\,2}_{2}\\
&\leq\,\int\limits_{\,X_{1}}\,\left<\,f,\, F_{1}\,(\,x_{\,1}\,)\,\right>_{1}\,\left<\,G_{1}\,(\,x_{\,1}\,),\, f\,\right>_{1}\,d\mu_{\,1}\,\int\limits_{\,X_{2}}\,\left<\,g,\, F_{2}\,(\,x_{\,2}\,)\,\right>_{2}\,\left<\,G_{2}\,(\,x_{\,2}\,),\, g\,\right>_{2}\,d\mu_{\,2}\\
&\,\leq\,B\,\left\|\,f\,\right\|^{\,2}_{1}\,\left\|\,g\,\right\|^{\,2}_{2}.
\end{align*} 
Since \,$f$\, and \,$g$\, are non-zero and therefore 
\[\int\limits_{\,X_{1}}\,\left<\,f,\, F_{1}\,(\,x_{\,1}\,)\,\right>_{1}\,\left<\,G_{1}\,(\,x_{\,1}\,),\, f\,\right>_{1}\,d\mu_{\,1}\,,\, \int\limits_{\,X_{2}}\,\left<\,g,\, F_{2}\,(\,x_{\,2}\,)\,\right>_{2}\,\left<\,G_{2}\,(\,x_{\,2}\,),\, g\,\right>_{2}\,d\mu_{\,2}\] are non-zero.\,Thus from the above inequality we can write
\begin{align*}
&\dfrac{A\,\left\|\,g\,\right\|^{\,2}_{2}}{\int\limits_{\,X_{2}}\,\left<\,g,\, F_{2}\,(\,x_{\,2}\,)\,\right>_{2}\,\left<\,G_{2}\,(\,x_{\,2}\,),\, g\,\right>_{2}\,d\mu_{\,2}}\,\left\|\,f\,\right\|^{\,2}_{1}\\
&\leq\,\int\limits_{\,X_{1}}\,\left<\,f,\, F_{1}\,(\,x_{\,1}\,)\,\right>_{1}\,\left<\,G_{1}\,(\,x_{\,1}\,),\, f\,\right>_{1}\,d\mu_{\,1}\\
&\leq\,\dfrac{B\,\left\|\,g\,\right\|^{\,2}_{2}}{\int\limits_{\,X_{2}}\,\left<\,g,\, F_{2}\,(\,x_{\,2}\,)\,\right>_{2}\,\left<\,G_{2}\,(\,x_{\,2}\,),\, g\,\right>_{2}\,d\mu_{\,2}}\,\left\|\,f\,\right\|^{\,2}_{1}.
\end{align*}
Thus, for each \,$f \,\in\, H_{1} \,-\, \{\,\theta\,\}$, we have
\begin{align*}
A_{\,1}\,\left\|\,f\,\right\|^{\,2}_{1} &\,\leq\, \int\limits_{\,X_{1}}\,\left<\,f,\, F_{1}\,(\,x_{\,1}\,)\,\right>_{1}\,\left<\,G_{1}\,(\,x_{\,1}\,),\, f\,\right>_{1}\,d\mu_{\,1} \leq\,B_{\,1}\,\left\|\,f\,\right\|^{\,2}_{1},
\end{align*}
where 
\[A_{\,1} \,=\, \inf_{g \,\in\, H_{\,2},\, \|\,g\,\|_{2} \,=\, 1}\,\left\{\,\dfrac{A\,\left\|\,g\,\right\|^{\,2}_{2}}{\int\limits_{\,X_{2}}\,\left<\,g,\, F_{2}\,(\,x_{\,2}\,)\,\right>_{2}\,\left<\,G_{2}\,(\,x_{\,2}\,),\, g\,\right>_{2}\,d\mu_{\,2}}\,\right\},\]and
\[B_{\,1} \,=\, \sup_{g \,\in\, H_{\,2},\, \|\,g\,\|_{2} \,=\, 1}\,\left\{\,\dfrac{B\,\left\|\,g\,\right\|^{\,2}_{2}}{\int\limits_{\,X_{2}}\,\left<\,g,\, F_{2}\,(\,x_{\,2}\,)\,\right>_{2}\,\left<\,G_{2}\,(\,x_{\,2}\,),\, g\,\right>_{2}\,d\mu_{\,2}}\,\right\}.\] 
This shows that \,$\left(\,F_{1},\,G_{1}\,\right)$\, is a continuous biframe for \,$H_{1}$\, with respect to \,$\left(\,X_{1},\, \mu_{\,1}\,\right)$. Similarly, it can be shown that \,$\left(\,F_{2},\,G_{2}\,\right)$\, is a continuous biframe for \,$H_{2}$\, with respect to \,$\left(\,X_{2},\, \mu_{\,2}\,\right)$.

Conversely, suppose that \,$\left(\,F_{1},\,G_{1}\,\right)$\, is a continuous biframe for \,$H_{1}$\, with respect to \,$\left(\,X_{1},\, \mu_{\,1}\,\right)$\, having bounds \,$A,\, B$\, and \,$\left(\,F_{2},\,G_{2}\,\right)$\, is a continuous biframe for \,$H_{2}$\, with respect to \,$\left(\,X_{2},\, \mu_{\,2}\,\right)$\, having bounds \,$C,\,D$.\,By the assumption it is easy to very that \,$\left(\,\mathbf{F},\, \mathbf{G}\,\right) \,=\,  \left(\,F_{1} \,\otimes\, F_{2},\, G_{1} \,\otimes\, G_{2}\,\right)$ is weakly measurable on \,$H_{1} \,\otimes\, H_{2}$\, with respect to \,$(\,X,\, \mu\,)$.\,Now, for each \,$f \,\in\, H_{1} \,-\, \{\,\theta\,\}$\, and \,$g \,\in\, H_{2} \,-\, \{\,\theta\,\}$, we have
\begin{align*}
&A\,\left\|\,f\,\right\|^{\,2}_{1} \,\leq\, \int\limits_{\,X_{1}}\,\left<\,f,\, F_{1}\,(\,x_{\,1}\,)\,\right>_{1}\,\left<\,G_{1}\,(\,x_{\,1}\,),\, f\,\right>_{1}\,d\mu_{\,1} \leq\,B\,\left\|\,f\,\right\|^{\,2}_{1}\,,\\
&C\,\left\|\,g\,\right\|^{\,2}_{2} \,\leq\, \int\limits_{\,X_{2}}\,\left<\,g,\, F_{2}\,(\,x_{\,2}\,)\,\right>_{2}\,\left<\,G_{2}\,(\,x_{\,2}\,),\, g\,\right>_{1}\,d\mu_{\,2} \leq\,B\,\left\|\,g\,\right\|^{\,2}_{2}.
\end{align*} 
Multiplying the above two inequalities and using Fubini's theorem we get
\begin{align*}
&A\,C\,\left\|\,f \,\otimes\, g\,\right\|^{\,2}\\
& \,\leq\, \int\limits_{\,X}\,\left<\,f \,\otimes\, g,\, F_{1}\,(\,x_{\,1}\,) \,\otimes\, F_{2}\,(\,x_{\,2}\,)\,\right>\,\left<\, G_{1}\,(\,x_{\,1}\,) \,\otimes\, G_{2}\,(\,x_{\,2}\,),\, f \,\otimes\, g\,\right>\,d\mu \\
&\leq\,B\,D\,\left\|\,f \,\otimes\, g\,\right\|^{\,2}, 
\end{align*} 
for all \,$f \,\otimes\, g \,\in\, H_{1} \,\otimes\, H_{2}$.\,Thus, for each \,$f \,\otimes\, g \,\in\, H_{1} \,\otimes\, H_{2}$, we have
\begin{align*}
&A\,C\,\left\|\,f \,\otimes\, g\,\right\|^{\,2} \,\leq\, \int\limits_{\,X}\,\left<\,f \,\otimes\, g,\, \mathbb{F}\,(\,x\,)\,\right>\,\left<\,\mathbb{G}\,(\,x\,),\, f \,\otimes\, g\,\right>\,d\mu \leq\,B\,D\,\left\|\,f \,\otimes\, g\,\right\|^{\,2},
\end{align*}
Hence, \,$\left(\,\mathbf{F},\, \mathbf{G}\,\right) \,=\,  \left(\,F_{1} \,\otimes\, F_{2},\, G_{1} \,\otimes\, G_{2}\,\right)$\, is a continuous biframe for \,$H_{1} \,\otimes\, H_{2}$\, with respect to \,$(\,X,\, \mu\,)$\, having bounds \,$A\,C$\, and \,$B\,D$.\,This completes the proof.  
\end{proof} 

\begin{example}
Let \,$\left\{\,e_{\,i}\,\right\}_{i \,=\, 1}^{\,\infty}$\, be an orthonormal basis for \,$H_{1}$\, and \,$\left(\,X_{1},\, \mu_{\,1}\,\right)$\, be a measure space with \,$\mu_{\,1}$\, is \,$\sigma$-finite.\,Then we can write \,$X_{1} \,=\, \bigcup_{i \,=\, 1}^{\,\infty}\,\Omega_{i}$, where \,$\left\{\,\Omega_{\,i}\,\right\}_{i \,=\, 1}^{\,\infty}$\, is a sequence of disjoint measurable subsets of \,$X_{1}$\, with \,$\mu\left(\,\Omega_{i}\,\right) \,<\, \infty$.\,Suppose
\begin{align*}
&\left\{\,f_{\,i}\,\right\}_{i \,=\, 1}^{\,\infty} \,=\, \left\{\,e_{\,1},\, e_{\,1},\, e_{\,1},\, e_{\,2},\, e_{\,3},\, \cdots\,\cdots\,\right\}\,,\\
&\left\{\,g_{\,i}\,\right\}_{i \,=\, 1}^{\,\infty} \,=\, \left\{\,0,\, e_{\,1},\, e_{\,1},\, e_{\,2},\, e_{\,3},\, \cdots\,\cdots\,\right\}.
\end{align*}
For each \,$x_{1} \,\in\, X_{1}$, we define the mappings \,$F_{1} \,:\, X_{1} \,\to\, H_{1}$\, by \,$F_{1}\,(\,x_{1}\,) \,=\, \dfrac{1}{\sqrt{\mu\left(\,\Omega_{i}\,\right)}}\,f_{\,i}$\, and \,$G_{1} \,:\, X_{1} \,\to\, H_{1}$\, by \,$G_{1}\,(\,x_{1}\,) \,=\, \dfrac{1}{\sqrt{\mu\left(\,\Omega_{i}\,\right)}}\,g_{\,i}$.
Then by Example \ref{3exm3.11}, \,$\left(\,F_{1},\, G_{1}\,\right)$\, is a continuous biframe for \,$H$\, with bounds \,$1$\, and \,$2$.

On the other hand, let \,$\left\{\,e^{\,\prime}_{\,j}\,\right\}_{j \,=\, 1}^{\,\infty}$\, be an orthonormal basis for \,$H_{2}$\, and \,$\left(\,X_{2},\, \mu_{\,2}\,\right)$\, be a measure space with \,$\mu_{\,2}$\, is \,$\sigma$-finite.\,Then \,$X_{2} \,=\, \bigcup_{i \,=\, 1}^{\,\infty}\,\Omega^{\,\prime}_{j}$, and \,$\left\{\,\Omega^{\,\prime}_{\,j}\,\right\}_{j \,=\, 1}^{\,\infty}$\, is a sequence of disjoint measurable subsets of \,$X_{2}$\, with \,$\mu\left(\,\Omega^{\,\prime}_{j}\,\right) \,<\, \infty$.\,Suppose 
\begin{align*}
&\left\{\,f^{\,\prime}_{\,j}\,\right\}_{i \,=\, 1}^{\,\infty} \,=\, \left\{\,5\,e_{\,1},\, 3\,e_{\,2},\, 2\,e_{\,3},\, 2\,e_{\,4},\, \cdots\,\cdots\,\right\}\,,\\
&\left\{\,g^{\,\prime}_{\,j}\,\right\}_{i \,=\, 1}^{\,\infty} \,=\, \left\{\,0,\, 0,\, 3\,e_{\,1},\, 0,\, 2\,e_{\,2},\, 2\,e_{\,3},\, \cdots\,\cdots\,\right\}.
\end{align*}
Now, we define \,$F_{\,2} \,:\, X_{2} \,\to\, H_{2}$\, by \,$F_{\,2}\,(\,x_{\,2}\,) \,=\, \dfrac{1}{\sqrt{\mu\left(\,\Omega^{\,\prime}_{j}\,\right)}}\,f^{\,\prime}_{\,j}$\, and  \,$G_{\,2} \,:\, X_{2} \,\to\, H_{2}$\, by \,$G_{\,2}\,(\,x_{\,2}\,) \,=\, \dfrac{1}{\sqrt{\mu\left(\,\Omega^{\,\prime}_{j}\,\right)}}\,g^{\,\prime}_{\,j}$.\,Now, for \,$f,\, g \,\in\, H_{2}$, we have
\begin{align*}
&\int\limits_{\,X_{2}}\,\left<\,f,\, F_{2}\,(\,x_{2}\,)\,\right>\,\left<\,G_{2}\,(\,x_{2}\,),\, f\,\right>\,d\mu_{2} \,=\,\sum\limits_{j \,=\, 1}^{\infty}\,\int\limits_{\,\Omega^{\,\prime}_{j}}\,\left <\,f,\, f^{\,\prime}_{\,j}\,\right >\,\left <\,g^{\,\prime}_{\,j},\, f\,\right >\,d\mu_{2}\\
&=\,\left <\,f,\, e_{\,1}\,\right >\,\left <\,e_{\,1},\, f\,\right > \,+\, 2\,\left <\,f,\, e_{\,1}\,\right >\,\left <\,e_{\,1},\, f\,\right > \,+\, 2\,\left <\,f,\, e_{\,2}\,\right >\,\left <\,e_{\,2},\, f\,\right >+\cdots\\
&=\,\left|\,\left <\,f,\, e_{\,1}\,\right >\,\right|^{\,2} \,+\, 2\,\left|\,\left <\,f,\, e_{\,1}\,\right >\,\right|^{\,2} \,+\, 2\,\left|\,\left <\,f,\, e_{\,2}\,\right >\,\right|^{\,2} \,+\, \cdots\\
& \,=\, \left|\,\left <\,f,\, e_{\,1}\,\right >\,\right|^{\,2} \,+\, 2\,\|\,f\,\|^{\,2}. 
\end{align*} 
Thus, \,$\left(\,F_{2},\, G_{2}\,\right)$\, is a continuous biframe for \,$H$\, with bounds \,$2$\, and \,$3$.\,Thus, by the Theorem \ref{th4.11}, \,$\left(\,\mathbf{F},\, \mathbf{G}\,\right) \,=\,  \left(\,F_{1} \,\otimes\, F_{2},\, G_{1} \,\otimes\, G_{2}\,\right)$\, is a continuous biframe for \,$H_{1} \,\otimes\, H_{2}$\, with respect to \,$(\,X,\, \mu\,)$\, having bounds \,$2$\, and \,$6$.
\end{example}

Let \,$(\,X,\, \mu\,) \,=\, \left(\,X_{1} \,\times\, X_{2},\, \mu_{\,1} \,\otimes\, \mu_{\,2}\,\right)$\, be the product of measure spaces with \,$\sigma$-finite positive measures \,$\mu_{\,1},\, \mu_{\,2}$.\,Let \,$L^{\,2}\left(\,X,\, \mu\,\right)$\, be the class of all measurable functions \,$\Psi \,:\, X \,\to\, H_{1} \,\otimes\, H_{2}$\, such that
\begin{align*}
&\int\limits_{\,X}\left\|\,\Psi\,(\,x\,)\,\right\|^{\,2}\,d\mu\\
&=\,\int\limits_{\,X_{1}}\,\left\|\,\varphi_{1}\,(\,x_{1}\,)\,\right\|_{\,1}\,d\mu_{1}\,\int\limits_{\,X_{2}}\,\left\|\,\varphi_{2}\,(\,x_{2}\,)\,\right\|_{\,2}\,d\mu_{2} \,<\, \infty,
\end{align*} 
for \,$\varphi_{1} \,\in\, L_{1}^{\,2}\left(\,X_{1},\, \mu_{1}\,\right)$,\, \,$\varphi_{2} \,\in\, L_{2}^{\,2}\left(\,X_{2},\, \mu_{2}\,\right)$, with the inner product 
\begin{align*}
\left<\,\Psi,\, \Phi\,\right>_{L^{\,2}} &\,=\, \int\limits_{\,X}\,\left<\,\Psi\,(\,x\,),\,\Phi\,(\,x\,)\,\right>\,d\mu \\
&=\,\int\limits_{\,X_{1}}\,\left<\,\varphi_{1}\,(\,x_{1}\,),\,\psi_{1}\,(\,x_{1}\,)\,\right>_{\,1}\,d\mu_{1}\,\int\limits_{\,X_{2}}\,\left<\,\varphi_{2}\,(\,x_{2}\,),\,\psi_{2}\,(\,x_{2}\,)\,\right>_{\,2}\,d\mu_{2}\\
&=\,\left<\,\varphi_{1},\, \psi_{1}\,\right>_{L_{1}^{\,2}}\,\left<\,\varphi_{2},\, \psi_{2}\,\right>_{L_{2}^{\,2}},
\end{align*}
where \,$\Psi \,=\, \varphi_{1} \,\otimes\, \varphi_{2}\,,\, \Phi \,=\, \psi_{1} \,\otimes\, \psi_{2} \,\in\, L^{\,2}\left(\,X,\, \mu\,\right)$, for \,$\varphi_{1},\,\psi_{1}  \,\in\, L_{1}^{\,2}\left(\,X_{1},\, \mu_{1}\,\right)$\, and \,$\varphi_{2},\,\psi_{2}  \,\in\, L_{2}^{\,2}\left(\,X_{2},\, \mu_{2}\,\right)$.\,The space \,$L^{\,2}\left(\,X,\, \mu\,\right)$\, is complete with respect to the above inner product.\,Therefore, it is an Hilbert space.

\begin{definition}
Let \,$(\,\mathbf{F},\, \mathbf{G}\,) \,=\, \,\left(\,\mathbb{F} \,:\, X \,\to\, H_{1} \,\otimes\, H_{2},\, \,\mathbb{G} \,:\, X \,\to\, H_{1} \,\otimes\, H_{2}\,\right)$\, be a continuous biframe for \,$H_{1} \,\otimes\, H_{2}$\, with respect to \,$(\,X,\, \mu\,)$.\,Then the operator \,$S_{\mathbb{F} \,\otimes\, \mathbb{G}} \,:\, H_{1} \,\otimes\, H_{2} \,\to\, H_{1} \,\otimes\, H_{2}$\, is given by
\[S_{\mathbb{F} \,\otimes\, \mathbb{G}}\,(\,f \,\otimes\, g\,) \,=\, \int\limits_{\,X}\,\left<\,f \,\otimes\, g,\, \mathbb{F}\,(\,x\,)\,\right>\,\mathbb{G}\,(\,x\,)\,d\mu\]is called continuous biframe operator.
\end{definition}

\begin{theorem}
Let \,$\left(\,\mathbf{F},\, \mathbf{G}\,\right) \,=\,  \left(\,F_{1} \,\otimes\, F_{2},\, G_{1} \,\otimes\, G_{2}\,\right) \,=\, \mathbb{F},\, \mathbb{G} \,:\, X \,\to\, H_{1} \,\otimes\, H_{2}$\, is a continuous biframe for \,$H_{1} \,\otimes\, H_{2}$\, with respect to \,$(\,X,\, \mu\,)$\, having continuous biframe operator \,$S_{\mathbb{F} \,\otimes\, \mathbb{G}}$.\,Then \,$S_{\mathbb{F} \,\otimes\, \mathbb{G}} \,=\, S_{F_{\,1},\, G_{1}}  \,\otimes\, S_{F_{\,2},\, G_{2}}$, where \,$S_{F_{\,1},\, G_{1}}$\, and \,$S_{F_{\,2},\, G_{2}} $\, are continuous biframe operators of \,$\left(\,F_{1},\, G_{1}\,\right)$\, and \,$\left(\,F_{2},\, G_{2}\,\right)$, respectively.   
\end{theorem}

\begin{proof}
For each \,$f \,\otimes\, g \,\in\, H_{1} \,\otimes\, H_{2}$, we have
\begin{align*}
&S_{\mathbb{F} \,\otimes\, \mathbb{G}}\,(\,f \,\otimes\, g\,)\\
 &= \int\limits_{X}\left<\,f \otimes g,\, F_{1}(\,x_{1}\,) \otimes F_{2}(\,x_{2}\,)\,\right>\,\left(\,G_{1}(\,x_{1}\,) \otimes G_{2}\,(\,x_{2}\,)\right)d\mu \\
&=\left(\,\int\limits_{\,X_{1}}\,\left<\,f,\, F_{1}\,(\,x_{\,1}\,) \,\right>_{1}\,G_{1}\,(\,x_{\,1}\,)\,d\mu_{\,1}\,\right)\,\otimes\left(\,\int\limits_{\,X_{2}}\,\left<\,g,\, F_{2}\,(\,x_{\,2}\,) \,\right>_{2}\,G_{2}\,(\,x_{\,2}\,)\,d\mu_{\,2}\,\right)\\
&=\, S_{F_{\,1},\, G_{1}}\,f  \,\otimes\, S_{F_{\,2},\, G_{2}}\,g \,=\, \left(\,S_{F_{\,1},\, G_{1}}  \,\otimes\, S_{F_{\,2},\, G_{2}}\,\right)\,(\,f \,\otimes\, g\,).
\end{align*} 
Thus, \,$S_{\mathbb{F} \,\otimes\, \mathbb{G}} \,=\, S_{F_{\,1},\, G_{1}}  \,\otimes\, S_{F_{\,2},\, G_{2}}$. 
\end{proof}

\begin{lemma}
Let \,$\left(\,F_{1},\, G_{1}\,\right)$\, be a continuous biframe for \,$H_{1}$\, with respect to \,$\left(\,X_{1},\, \mu_{1}\,\right)$\, having bounds \,$A,\, B$\, and \,$\left(\,F_{2},\, G_{2}\,\right)$\, be a continuous biframe for \,$H_{2}$\, with respect to \,$\left(\,X_{2},\, \mu_{2}\,\right)$\, having bounds \,$C,\,D$.\,Then \,$A\,C\,I_{H_{1} \,\otimes\, H_{2}} \,\leq\, S_{\mathbb{F} \,\otimes\, \mathbb{G}} \,\leq\, B\,D\,I_{H_{1} \,\otimes\, H_{2}}$, where \,$\,I_{H_{1} \,\otimes\, H_{2}}$\, is the identity operator on \,$H_{1} \,\otimes\, H_{2}$\, and \,$S_{F_{1},\, G_{1}},\, \,S_{F_{2},\, G_{2}}$\, are continuous biframe operators of \,$\left(\,F_{1},\, G_{1}\,\right),\, \left(\,F_{2},\, G_{2}\,\right)$, respectively.
\end{lemma}

\begin{proof}
Since \,$S_{F_{1},\, G_{1}}$\, and \,$S_{F_{2},\, G_{2}}$\, are continuous biframe operators, we have 
\[A\,I_{H_{1}} \,\leq\, S_{F_{1},\, G_{1}} \,\leq\, B\,I_{H_{1}},\; \;C\,I_{H_{2}} \,\leq\, S_{F_{2},\, G_{2}} \,\leq\, D\,I_{H_{2}},\]
where \,$I_{H_{1}}$\, and \,$I_{H_{2}}$\, are the identity operators on \,$H_{1}$\, and \,$H_{2}$, respectively.\,Taking tensor product on the above two inequalities, we get
\begin{align*}
&A\,C \left(\,I_{H_{1}} \,\otimes\, I_{H_{2}}\,\right) \,\leq\,  \left(\,S_{F_{1},\, G_{1}} \,\otimes\, S_{F_{1},\, G_{1}}\,\right) \,\leq\, B\,D\,\left(\,I_{H_{1}} \,\otimes\, I_{H_{2}}\,\right)\\
&\Rightarrow\,A\,C\,I_{H_{1} \,\otimes\, H_{2}} \,\leq\, S_{\mathbb{F} \,\otimes\, \mathbb{G}} \,\leq\, B\,D\,I_{H_{1} \,\otimes\, H_{2}}
. 
\end{align*}
This completes the proof.
\end{proof}

\begin{theorem}
Let \,$\left(\,F_{1},\, G_{1}\,\right)$\, be a continuous biframe for \,$H_{1}$\, with respect to \,$\left(\,X_{1},\, \mu_{1}\,\right)$\, and \,$\left(\,F_{2},\, G_{2}\,\right)$\, be a continuous biframe for \,$H_{2}$\, with respect to \,$\left(\,X_{2},\, \mu_{2}\,\right)$. Then \,$\Delta \,=\, \left(\,\left(\,T_{1} \,\otimes\, T_{2}\,\right)\left(\,F_{1}\,\otimes\, F_{2}\,\right),\, \left(\,T_{1} \,\otimes\, T_{2}\,\right)\left(\,G_{1}\,\otimes\, G_{2}\,\right)\,\right)$\, is a continuous biframe for \,$H_{1} \,\otimes\, H_{2}$\, with respect to \,$\left(\,X,\, \mu\,\right)$\, if and only if \,$T_{1} \,\otimes\, T_{2}$\, is an invertible bounded linear operator on \,$H_{1} \,\otimes\, H_{2}$.  
\end{theorem}

\begin{proof}
Let \,$\left(\,F_{1},\, G_{1}\,\right)$\, be a continuous biframe for \,$H_{1}$\, with bounds \,$A,\, B$\, and \,$\left(\,F_{2},\, G_{2}\,\right)$\, be a continuous biframe for \,$H_{2}$\, with bounds \,$C,\, D$.\,First we suppose that \,$T_{1} \,\otimes\, T_{2}$\, is an invertible bounded linear operator on \,$H_{1} \,\otimes\, H_{2}$.\,Then by Theorem \ref{th1.1}, \,$T_{1}$\, and \,$T_{2}$\, are invertible bounded linear operators on \,$H_{1}$\, and \,$H_{2}$, respectively.\,Now, by Theorem \ref{3.thm3.39}, \,$\left(\,T_{1}F_{1},\, T_{1}G_{1}\,\right)$\, is a continuous biframe for \,$H_{1}$\, with bounds \,$A\,\left\|\,T_{1}^{\,-\, 1}\,\right\|^{\,-\, 2},\, B\,\left\|\,T_{1}\,\right\|^{\,2}$\, and \,$\left(\,T_{2}F_{2},\, T_{2}G_{2}\,\right)$\, is a continuous biframe for \,$H_{2}$\, with bounds \,$C\,\left\|\,T_{2}^{\,-\, 1}\,\right\|^{\,-\, 2},\, D\,\left\|\,T_{2}\,\right\|^{\,2}$.\,Hence, by Theorem \ref{th4.11}, the pair
\begin{align*}
&\left(\,T_{1}F_{1} \,\otimes\, T_{2}F_{2}\,,\, T_{1}G_{1} \,\otimes\, T_{2}G_{2}\,\right)\\
& \,=\, \left(\,\left(\,T_{1} \,\otimes\, T_{2}\,\right)\left(\,F_{1}\,\otimes\, F_{2}\,\right),\, \left(\,T_{1} \,\otimes\, T_{2}\,\right)\left(\,G_{1}\,\otimes\, G_{2}\,\right)\,\right)
\end{align*} 
is a continuous biframe for \,$H_{1} \,\otimes\, H_{2}$\, with bounds 
\begin{align*}
\dfrac{A\,C}{\left\|\,T_{1}^{\,-\, 1}\,\right\|^{\,2}\,\left\|\,T_{2}^{\,-\, 1}\,\right\|^{\,2}} \,=\, \dfrac{A\,C}{\left\|\,T_{1}^{\,-\, 1} \,\otimes\, T_{2}^{\,-\, 1}\,\right\|^{\,2}}& \,=\, \dfrac{A\,C}{\left\|\,\left(\,T_{1} \,\otimes\, T_{2}\,\right)^{\,-\, 1}\,\right\|^{\,2}} \\
&\,=\, A\,C\,\left\|\,\left(\,T_{1} \,\otimes\, T_{2}\,\right)^{\,-\, 1}\,\right\|^{\,-\, 2}
\end{align*} 
and \,$B\,D\,\left\|\,T_{1}\,\right\|^{\,2}\,\left\|\,T_{2}\,\right\|^{\,2} \,=\, B\,D\,\left\|\,T_{1} \,\otimes\, T_{2}\,\right\|^{\, 2}$. 

Conversely, suppose that \,$\Delta$\, is a continuous biframe for \,$H_{1} \,\otimes\, H_{2}$\, with respect to \,$\left(\,X,\, \mu\,\right)$.\,Let \,$S_{F_{1},\, G_{1}}$\, and \,$S_{F_{2},\, G_{2}}$\, be the continuous biframe operators of \,$\left(\,F_{1},\, G_{1}\,\right)$\, and \,$\left(\,F_{2},\, G_{2}\,\right)$, respectively.\,By Theorem \ref{th4.11}, \,$\left(\,T_{1}F_{1},\, T_{1}G_{1}\,\right)$\, and \,$\left(\,T_{2}F_{2},\, T_{2}G_{2}\,\right)$\, are continuous biframes for \,$H_{1}$\, and \,$H_{2}$, respectively.\,Now, for \,$f \,\in\, H_{1}$, we have 
\begin{align*}
&\int\limits_{\,X_{1}}\,\left<\,f,\, T_{1}\,F_{1}\,(\,x_{\,1}\,) \,\right>_{1}\,T_{1}\,G_{1}\,(\,x_{\,1}\,)\,d\mu_{\,1}\\
& \,=\, T_{1}\,\left(\,\int\limits_{\,X_{1}}\,\left<\,T^{\,\ast}_{1}\,f,\, F_{1}\,(\,x_{\,1}\,) \,\right>_{1}\,G_{1}\,(\,x_{\,1}\,)\,d\mu_{\,1}\,\right) \,=\, T_{1}\,S_{F_{1},\, G_{1}}\,T^{\,\ast}_{1}\,f.
\end{align*} 
This shows that \,$T_{1}\,S_{F_{1},\, G_{1}}\,T^{\,\ast}_{1}$\, is the corresponding continuous biframe operator of \,$\left(\,T_{1}F_{1},\, T_{1}G_{1}\,\right)$.\,Thus,  \,$T_{1}\,S_{F_{1},\, G_{1}}\,T^{\,\ast}_{1}$\, is invertible on \,$H_{1}$\, and hence \,$T_{1}$\, is invertible on \,$H_{1}$.\,Similarly, it can be shown that \,$T_{2}\,S_{F_{2},\, G_{2}}\,T^{\,\ast}_{2}$\, is the corresponding continuous biframe operator of \,$\left(\,T_{2}F_{2},\, T_{2}G_{2}\,\right)$\, and hence \,$T_{2}$\, is invertible on \,$H_{2}$.\,Therefore, \,$T_{1} \,\otimes\, T_{2}$\, is an invertible bounded linear operator on \,$H_{1} \,\otimes\, H_{2}$.\,This completes the proof.      
\end{proof}

Now, we would complete this section by introducing the idea of a continuous biframe Bessel multiplier in \,$H_{1} \,\otimes\, H_{2}$.   

\begin{definition}
Let \,$(\,\mathbb{F},\, \mathbb{F}\,)$\, and \,$(\,\mathbb{G},\, \mathbb{G}\,)$\, be continuous biframe Bessel mappings in \,$H_{1} \,\otimes\, H_{2}$\, with respect to \,$\left(\,X,\, \mu\,\right)$\, having bounds \,$B_{1}$\, and \,$B_{2}$, respectively and \,$m \,:\, X \,\to\, \mathbb{C}$\, be a measurable function.\,The operator \,$\mathcal{M}_{m,\, \mathbb{F},\, \mathbb{G}} \,:\, H_{1} \,\otimes\, H_{2} \,\to\, H_{1} \,\otimes\, H_{2}$\, defined by
\begin{align}
&\mathcal{M}_{m,\, \mathbb{F},\, \mathbb{G}}\,(\,f \,\otimes\, g\,)\nonumber\\
& \,=\, \int\limits_{\,X}\,m\,(\,x\,)\,\left<\,f \,\otimes\, g,\, \mathbb{F}\,(\,x\,)\,\right>\,\mathbb{G}\,(\,x\,)\,d\mu\,,\label{eqn4.23}
\end{align} 
for all \,$f \,\otimes\, g \,\in\, H_{1} \,\otimes\, H_{2}$\,, is called continuous biframe Bessel multiplier of \,$\mathbb{F}$\, and \,$\mathbb{G}$\, with respect to \,$m$. 
\end{definition}

\begin{note}
Let \,$F_{1}\,,\, G_{1} \,:\, X_{1} \,\to\, H_{1}$\, be continuous biframe Bessel mappings in \,$H_{1}$\, with respect to \,$\left(\,X_{1},\, \mu_{1}\,\right)$\, and \,$F_{2}\,,\, G_{2} \,:\, X_{2} \,\to\, H_{2}$\, be continuous biframe Bessel mappings in \,$H_{2}$\, with respect to \,$\left(\,X_{2},\, \mu_{2}\,\right)$\, and \,$m_{1} \,:\, X_{1} \,\to\, \mathbb{C}$, \,$m_{2} \,:\, X_{2} \,\to\, \mathbb{C}$\, be two measurable functions.\,Suppose \,$M_{m_{1},\, F_{1},\, G_{1}} \,:\, H_{1} \,\to\, H_{1}$\, be the continuous biframe Bessel multiplier of \,$F_{1}$\, and \,$G_{1}$\, with respect to \,$m_{1}$\, and \,$M_{m_{2},\, F_{2},\, G_{2}} \,:\, H_{2} \,\to\, H_{2}$\, be the continuous biframe Bessel multiplier of \,$F_{2}$\, and \,$G_{2}$\, with respect to \,$m_{2}$.\,Now, by Theorem \ref{th4.11}, \,$\mathcal{F} \,=\, F_{1} \,\otimes\, F_{2}$, \,$\mathcal{G} \,=\, G_{1} \,\otimes\, G_{2}  \,:\, X \,\to\, H_{1} \,\otimes\, H_{2}$\, are continuous biframe Bessel mappings in \,$H_{1} \,\otimes\, H_{2}$\, with respect to \,$(\,X,\, \mu\,)$.\,Also, \,$m\,(\,x\,) \,=\, m_{1}\,(\,x_{1}\,)\,m_{2}\,(\,x_{2}\,)$\, is a measurable function.\,From (\ref{eqn4.23}), for each \,$f \,\otimes\, g \,\in\, H_{1} \,\otimes\, H_{2}$, we can write
\begin{align*}
&\mathcal{M}_{m,\, \mathbb{F},\, \mathbb{G}}\,(\,f \,\otimes\, g\,)\\
& \,=\, \int\limits_{\,X}\,m\,(\,x\,)\,\left<\,f \,\otimes\, g,\, \mathbb{F}\,(\,x\,)\,\right>\,\mathbb{G}\,(\,x\,)\,d\mu\\
&= \int\limits_{\,X_{1}}\,m_{1}(\,x_{1}\,)\left<f,\, F_{1}\,(\,x_{\,1}\,)\right>_{1}G_{1}(\,x_{\,1}\,)\,d\mu_{\,1}\,\otimes\int\limits_{X_{2}}\,m_{2}(\,x_{2}\,)\left<g,\, F_{2}\,(\,x_{\,2}\,)\right>_{2}G_{2}(\,x_{\,2}\,)\,d\mu_{\,2}\\
&=\, M_{m_{1},\, F_{1},\, G_{1}}\,f \,\otimes\, M_{m_{2},\, F_{2},\, G_{2}}\,g \,=\, \left(\, M_{m_{1},\, F_{1},\, G_{1}} \,\otimes\, M_{m_{2},\, F_{2},\, G_{2}}\,\right)\,(\,f \,\otimes\, g\,).
\end{align*}
Thus, \,$\mathcal{M}_{m,\, \mathbb{F},\, \mathbb{G}} \,=\, M_{m_{1},\, F_{1},\, G_{1}} \,\otimes\, M_{m_{2},\, F_{2},\, G_{2}}$. 
\end{note}

\begin{remark}
Following the Theorem \ref{thm4.22}, the continuous biframe Bessel multiplier of \,$\mathbb{F}$\, and \,$\mathbb{G}$\, with respect to \,$m$\, is well defined and bounded. 
\end{remark}

\end{document}